\magnification 1200
\input amssym.def
\input amssym.tex

\parindent = 40 pt
\parskip = 12 pt
\font \heading = cmbx10 at 12 true pt
 at 22 true pt
\font \medheading =cmbx7 at 15 true pt
 at 7 true pt
\def \R{{\bf R}}
\def \C{{\bf C}}
\def \Q{{\bf Q}}

\centerline {\medheading  A Constructive Elementary Method for Local}
\centerline {\medheading  Resolution of Singularities }
\rm
\line{}
\line{}
\centerline{\heading Michael Greenblatt}
\line{}
\centerline{December 30, 2014}
\baselineskip = 12 pt
\font \heading = cmbx10 at 14 true pt
\line{}
\line{}
\noindent{\heading 1. Introduction and statements of results}

\vfootnote{}{This research was supported in part by NSF grants DMS-0919713 and DMS-1001070} In this paper we
considerably simplify the local  resolution of singularities techniques of [G1]-[G3] and extend it to functions with convergent power 
series over an arbitrary local field $K$ of characteristic zero. By the classification of such fields, this means $K = \R, \C$, the $p$-adics $\Q_p$, or a finite field extension of $\Q_p$. The algorithm of this paper is an entirely elementary
and self-contained classical analysis argument, using only the implicit function theorem and elementary facts about power series and Newton polyhedra, and 
the methods are quite different from traditional resolution of singularities methods. In a
separate paper [G4], the methods of this paper (but for the most part not the resolution of singularities theorems themselves) are used to prove new results concerning exponential sums, oscillatory integrals, and related matters. 

The blowups used in this paper will all be traditional blowups along $n-2$ dimensional coordinate hyperplanes in appropriate coordinate systems. The reason for this is that such blowups correspond to 
the simplest permissible linear transformations on the Newton polyhedron of the function whose zero set is being resolved, and 
linear transformations on such Newton polyhedra play a major role in our arguments. Thus such blowups form natural building blocks
for the more complicated linear transformations we will need. They are also relatively easily made compatible with partitions of
unity, as we will see.  

\noindent The coordinate changes
used in the resolution process will be of the following types.

\noindent {\bf Definition 1.1.} An {\it invertible affine linear map} is a map $Ax + b$ on $K^n$, where $A$ is linear and invertible
and where $b \in K^n$.

\noindent {\bf Definition 1.2.} A {\it blowup} is a map $b:K^n \rightarrow K^n$ of the form $b(x) = 
(b_1(x),...,b_n(x))$ where there exist $j,k$ such that $b_i(x) = x_i$ for all $i \neq j$ and $b_j(x) = x_jx_k$.

\noindent {\bf Definition 1.3.} A {\it quasitranslation} is a function on a neighborhood of the origin in $K^n$ of the form
$$g(x) =  (x_1,...,x_{i-1},x_i + a(x_1,...,x_{i-1},x_{i+1},...,x_n),x_{i+1},...,x_n)$$
Here $a(x_1,...,x_{i-1},x_{i+1},...,x_n)$ has a convergent power series on a neighborhood of the origin in $K^{n-1}$
with $a(0) = 0$.

\noindent Our local resolution of singularities theorem is as follows. 

\noindent {\bf Theorem 1.1.}  Let $K$ be a local field of characteristic zero. Let $f(x)  = \sum_{\alpha}f_{\alpha} x^{\alpha}$ be a function with 
a convergent power series on an open $U_0 \subset K^n$ containing the origin such 
that $f(0) = 0$ with at least one $f_{\alpha}$ being nonzero. Let $K_0$ denote $K - \{0\}$. There is a bounded open $U$  with
$0 \in U \subset U_0$ such that if 
$\phi(x)$ is a smooth function supported in $U$ then outside a set of measure zero $\phi(x)$ can be written as $\sum_{i=1}^N\phi_i(x)$, each $\phi_i(x)$ supported in $U$, where for each $i$ there 
is a bounded open set $W_i \subset K^n$ and a function $\alpha_i:W_i \rightarrow U$ satisfying the following.

\noindent {\bf a)} Each $\alpha_i$ is a
finite composition of invertible affine linear maps, blowups, and quasitranslations such that the restriction of $\alpha_i$ to $W_i \cap K_0^n$ is an analytic diffeomorphism onto its image.

\noindent {\bf b)} Let $V_i = \{x \in W_i: \phi_i 
\circ \alpha_i(x) \neq 0\}$. Then $V_i$ is an open subset of $W_i \cap K_0^n$ with $cl(V_i) \subset W_i$. One can change $\phi_i \circ \alpha_i(x)$ 
on $W_i - K_0^n$ in such a way that the resulting function is smooth. Furthermore, $\alpha_i$ is an analytic diffeomorphism from $V_i$ to $\{x \in U: \phi_i(x) \neq 0\}$.

\noindent {\bf c)} On $W_i$, $f \circ \alpha_i(x)$  and the 
Jacobian determinant of $\alpha_i(x)$ can both be written in the form $c_i(x)m_i(x)$, where $c_i(x)$ is nonvanishing and $m_i(x)$ is a monomial.

\noindent Note that in the case where $K$ is a $p$-adic field, after adjusting on $W_i - K_0^n$, the functions $\phi_i \circ \alpha_i(x)$ are locally constant. 

Theorem 1.1 directly implies a more traditionally stated local resolution of singularities theorem involving open sets covering a an open set containing the origin:

\noindent {\bf Theorem 1.2.} Let $f(x)$ and $U_0$ be as in Theorem 1.1. There are  bounded open sets $V_i \subset K^n$, $i = 1,...,N$ such that for each $i$ there is a function $\alpha_i:V_i \rightarrow U_0$ satisfying the following.

\noindent {\bf a)} Each $\alpha_i$ is a finite composition of invertible affine linear maps, blowups, and quasitranslations which is
an analytic diffeomorphism from $V_i$ onto its image. 

\noindent {\bf b)} Up to a set of measure zero, $\cup_{i=1}^N \alpha_i(V_i)$ contains an open set containing the origin.

\noindent {\bf c)} Each $\alpha_i$ extends to an open $W_i$ containing $cl(V_i)$ such that on $W_i$, both $f \circ \alpha_i(x)$  and the Jacobian determinant of $\alpha_i(x)$ can  be written in the form $c_i(x)m_i(x)$, where $c_i(x)$ is nonvanishing and $m_i(x)$ is a monomial.

Theorem 1.2 follows immediately from Theorem 1.1 by letting $\phi(x)$ be any function equal to $1$ on a neighborhood of the
origin and then letting $V_i$ be as in Theorem 1.1. 

Theorems 1.1 and 1.2, being local results, are of course weaker than Hironaka's famous work [H1]-[H2] and their analytic counterparts. It should be pointed out there have
been other elementary proofs of resolution of singularities theorems over $\R$ and $\C$, the most well-known being
 [BM1] [Pa1]-[Pa2] [Su]. We also refer to the recent paper [CoGrPr] for further elementary resolution of singularities
results, modeled on [J] [Pa2] [Su]. It should also be mentioned that there have also been a number of significant simplifications and extensions of Hironaka's proof in recent years, such as [BM2] [EH] [EV1]-[EV2] [K] [W]. While these papers are in the language of modern algebraic geometry and thus are not classical analysis arguments, they are quite accessible and give quite general results of an algebraic nature.

This paper will proceed as follows. (The following is an oversimplication, but will give the reader an idea
of the overall argument). Let $f(x)$ have a convergent power series
on a neighborhood of the origin, not identically zero, such that $f(0) = 0$. In section 2, we will subdivide a small neighborhood of the origin into a finite 
collection of "wedges" $A_{ij}$. To each wedge $A_{ij}$ there corresponds a unique vertex or compact face of the Newton polyhedron
 of $f(x)$ of dimension $i$  such that the terms of $f(x)$'s Taylor expansion whose exponents are on that face dominate the overall
Taylor expansion of $f(x)$ in a certain way that is made precise in Theorem 2.1.

In section 3, we will take a small cube $F$ centered at the origin and write $F = \cup_{k = 1}^q F_k$ such that
after a sequence of blowups, $F_k$ becomes another cube $E_k$. These blowups are done such that (Lemma 3.1) the functions
$x^v$ for $v$  a vertex of the Newton polyhedron of $f$ become, in the blown up coordinates, of the form $x^w$ where the 
$w$ are now ordered lexicographically. This has the effect that the intersection of each $E_k$ with each  blown-up
$A_{ij}$ is either empty or is a set $A_{ijk}'$ comparable in a certain sense to a set of the form $C \times D$, where 
$C$ is a cube centered at the origin of  dimension $n - i$ and where $D$ is an open set of dimension of $i$ which does not 
intersect the coordinate hyperplanes. If $i = 0$ then a nonempty $A_{ijk}'$ is just comparable to a cube $C$ of dimension $n$. The 
exact statements (and what "comparable" means here) are given in  Lemmas 3.2 and 3.2'.
We will need
$A_{ijk}'$ to be of this type of regular shape for our induction to proceed.

We will also see (Lemma 3.5) that the above operations can be made compatible with partitions of unity; if $\phi(x)$ is smooth and supported on a small enough neighborhood of the origin, one can write $\phi(x) = \sum_{i=1}^N \phi_i(x)$ such that each
 $\phi_i(x)$ transforms under one of the above sequences of blowups  to a smooth function that is approximately supported on one of the $A_{ijk}'$. 

In section 4, we will prove Theorem 1.1. The beginning of the proof  is similar to how Hironaka-based proofs work. We induct on the dimension $n$, and for a given $n$ we induct on the order $m$ of the zero of $f(x)$ at the origin. We rotate coordinates so that
${\partial^m f \over \partial x_n^m}(0) \neq 0$ and ${\partial^l f \over \partial x_n^l}(0)  = 0$ for $l < m$. We then do a quasitranslation after which ${\partial^{m-1} f \over \partial x_n^{m-1}}(x_1,...,x_{n-1},0) = 0$ for all $x_1,...,x_{n-1}$, and 
then use the induction hypothesis on the dimension $n$, applied to each ${\partial^{l} f \over \partial x_n^{l}}(x_1,...,x_{n-1},0)$
with $l < m$, and reduce to the case where $f(x)$ is of the form
$$f(x)  = a(x_1,...,x_n) x_n^m + \sum_{p < m-1}a_{p}(x_1,...,x_{n-1})m_{p}(x_1,...,x_{n-1})x_n^p \eqno (1.1)$$
Here each $m_p(x_1,...,x_{n-1})$ is a monomial and $ a(x_1,...,x_n) $ and each $a_{p}(x_1,...,x_{n-1})$ are nonvanishing 
near the origin.

The argument now diverges from typical Hironaka-based arguments. With $f(x)$ now of the form $(1.1)$,
we divide a neighborhood of the origin into the above-mentioned $F_k$ and perform the associated blowups converting $F_k$
to $E_k$. Using the results of Section 3, we will see that because of the way the 
$A_{ij}$  are defined, the fact that ${\partial^m f \over \partial x_n^m} \neq 0$ and the fact that there are no terms in the Taylor expansion of $f(x)$ for which $x_n$ appears to degree $m-1$ imply that on each $A_{ijk}'$ we can
factor out a monomial after which we are left with a function $g(x)$ whose zeroes have order at most $m - 1$. 

Then induction on $m$ will give the result; $g(x)$ will satisfy Theorem 1.1 and our methods will be such 
 that the monomial being  factored out also will blow up to a function of the desired form $c(x)m(x)$ under the coordinate changes on $g(x)$, where $m(x)$ is a monomial and $c(x)$ is nonvanishing. Thus we end out with a sequence of
blowups after which $f(x)$ is also of this form. An inspection of the argument then reveals that the Jacobian of the composition of all coordinate
changes is also of the desired form  $c(x)m(x)$. Furthermore, the overall argument will be performed so as to be compatible with partitions of unity as required by Theorem 1.1.

 It should be pointed out that while
applying the induction hypothesis on the dimension on  coefficient functions such as in $(1.1)$ frequently occurs in resolution
of singularities algorithms that are related to Hironaka's work, the divisions into domains with respect to the Newton
polyhedron in sections 2 and the subsequent arguments of section 3 are different and are important parts of the proof of Theorem 1.1. The results of
sections 2 and 3 are also used in the analysis of oscillatory integrals and related matters in [G3] and [G4].  Such divisions have
antecedents, such as in the use of toric resolution of singularities
in the analysis of oscillatory integrals such as in [V], as well as the various Phong-Stein-Sturm papers on oscillatory integral 
operators in two dimensions such as [PS] and [PSSt]. 

One aspect of our arguments worth mentioning is that unlike in
traditional resolution of singularities methods, we do not require the centers of the blowups to be contained in the singular
set or the zero set of the function being resolved. The classical analysis applications of this paper so far have not required it, and not 
making this stipulation gives us more flexibility in our proofs.

\noindent {\heading 2. Subdivisions according to the Newton polyhedron.}

\noindent We start with some relevant definitions.

\noindent {\bf Definition 2.1.} Let $f(x)$ be a function such that $f(x)$ has a convergent power series expansion
$\sum_{\alpha} f_{\alpha}x^{\alpha}$ on a neighborhood of the origin in $K^n$.
For any $\alpha$ for which $f_{\alpha} \neq 0$, let $Q_{\alpha}$ be the orthant $\{t \in \R^n: 
t_i \geq \alpha_i$ for all $i\}$. Then the {\it Newton polyhedron} $N(f)$ of $f(x)$ is defined to be 
the convex hull of all $Q_{\alpha}$.  

A Newton polyhedron can contain faces of various dimensions in various configurations. 
These faces can be either compact or unbounded. In this paper, as in earlier work such as [G3] and [V], an 
important role is played by the following functions, defined for compact faces of the Newton polyhedron. A vertex
 is always considered to be a compact face of dimension zero.

\noindent {\bf Definition 2.2.} Suppose $F$ is a compact face of $N(f)$. Then
if $f(x) = \sum_{\alpha} f_{\alpha}x^{\alpha}$ denotes the Taylor expansion of $f$ like above, 
define $f_F(x) = \sum_{\alpha \in F} f_{\alpha}x^{\alpha}$.

\noindent {\bf Definition 2.3.} $v(f)$ denotes the set of vertices of $N(f)$.

We now describe how we subdivide a neighborhood of the origin according to the Newton polyhedron of a function.
Let $U$ be an open set in $K^n$ containing the origin and let $f:U \rightarrow K$ be a function with a (nontrivial) convergent power series on $U$ such that $f(0,...,0) = 0$.  As above, let $N(f)$ denote the Newton polyhedron of $f$, and let
$v(f)$ denote the set of vertices of $N(f)$. We list the compact faces of $N(f)$ of dimension $i$ as $\{F_{ij}\}_{j=1}^{k_i}$;
a vertex of $N(f)$ is considered to be a face of dimension zero. In the following, we will make use of constants $C_1, ..., C_n$
such that $C_1 > N$ and each $C_{i+1} > C_i^N$ for a large $N$ dictated by our arguments that will depend only on $N(f)$.

Letting $K_0 = K - \{0\}$ as before, we will write the cube $E = 
\{x \in  K_0^n: 0 < |x_l| < C_n^{-1}$ for all $l \}$ as the union, modulo a set of measure zero, of regions
$A_{ij}$. On a given $A_{ij}$, the terms $f_{\alpha}x^{\alpha}$ for $\alpha \in F_{ij}$ will dominate the Taylor expansion of $f(x)$ in a sense to be made precise by Theorem 2.1 and Lemma 2.2. 
We will start by defining $A_{n-1 j}$ successively in $j$. Once the $A_{ij}$ have been defined for a given $i$, we will 
proceed to define the $A_{i-1j}$ successively in $j$, and so on, until we have defined the $A_{0j}$. For a given $i$, it does
not matter in what order we list the faces $F_{ij}$ here as long as we have fixed some enumeration.

The definition of the $A_{ij}$ is as follows. Once all $A_{i'j'}$ have been defined for $i' > i$ and for $i' = i$, $j' < j$, 
we define $A_{ij}$ to be the set of points $x$ in the cube $E$ such that the following three conditions hold.

\noindent a) $\sup_{v \in v(f) \cap F_{ij}}|x^{v}| = \sup_{v \in v(f)}|x^v|$.

\noindent b)  $\inf_{v \in v(f) \cap F_{ij}}|x^{v}| \geq C_i^{-1}\sup_{v \in v(f)}|x^v|\,\,\,(= C_i^{-1}\sup_{v \in  v(f) \cap F_{ij}}|x^v|)$.

\noindent c) $x \notin A_{i'j'}$ for any $(i',j')$ with $i' > i$ or with $i' = i$ and $j' < j$.

Condition c) is omitted when defining $A_{n-1,1}$. Note that in the case where $i = 0$ condition b) is tautologically satisfied, so 
that the $x$ in $A_{0j}$ are those $x$ for which a specific $|x^v|$ is maximal and which were not selected to be in any 
previous $A_{ij}$. 

We next prove a lemma which will be important in the proof of the main theorem of this section, Theorem 2.1.

\noindent {\bf Lemma 2.0.}  Suppose $F_{ij}$ is a face of $N(f)$ and $F$ is either a face of $N(f)$ properly containing $F_{ij}$
or is $N(f)$ itself. Let $w$ be a point in the interior of $F$. Then there is a $\mu_0> 0$ depending on $w$ such that if the number $N$ in the
definition of the $C_i$ were chosen sufficiently large (also depending on $w$), then we have $|x^{w}| < C_{i+1}^{-\mu_0}\inf_{v \in v(f) 
\cap F_{ij}}|x^v|$ for all $x \in A_{ij}$.

\noindent {\bf Proof.} The 
argument is slightly different when $F$ is bounded from when it is unbounded, so we will consider the two cases separately, 
starting with the unbounded case.

\noindent {\bf Case 1. $F$ is unbounded.} 

When $F$ is unbounded there is some direction $c = (c_1,...,c_n)$ tangent to $F$ such that each component $c_i$ is nonnegative,
with at least one positive. Since $w$ is in the interior of $F$, we may replace $c$ by a small multiple of itself and say that 
$w - c \in F$. But since $w- c \in N(f)$, one can necessarily write
$$w - c= \sum_{v_k \in v(f)} \alpha_k v_k + d \eqno (2.1)$$
Here each $\alpha_k$ is nonnegative, $\sum_k \alpha_k = 1$, and all the components of $d$ are nonnegative.
Therefore
$$|x^w| = |x^{c+d}|\prod_{v_k \in v(f)}|x^{v_k}|^{\alpha_k} \eqno (2.2)$$
If $x \in A_{ij}$ then each factor $|x^{v_k}|$ is at most $C_i \inf_{v \in v(f) \cap F_{ij}}|x^v|$  by condition
b) of the definition of $A_{ij}$. As a result, since $\sum_k \alpha_k = 1$ equation $(2.2)$ implies 
$$|x^w| \leq C_i|x^{c+d}|\inf_{v \in v(f) \cap F_{ij}}|x^v|\eqno (2.3)$$
Since at least one component of $c+d$ is positive, if the diameter $C_n^{-1}$ of the original $E$ were chosen sufficiently small, for 
some $\mu_0 > 0$ we have 
$$C_i|x^{c+d}|< C_{i+1}^{-\mu_0}$$
This gives $|x^w| < C_{i+1}^{-\mu_0}\inf_{v \in v(f) \cap F_{ij}}|x^v|$  as needed. Thus we are done with the lemma in the
 case where $F$ is unbounded.

\noindent {\bf Case 2. $F$ is bounded.}

We now move onto the case where $F$ is a bounded, and thus compact face of $N(f)$. Let $(i',j')$ be such that $F = F_{i'j'}$.
Since $w$ is in the interior of $F$, one can write 
$$w = \sum_{v_k \in v(f) \cap F_{i'j'}}\alpha_kv_k \eqno (2.4)$$
Here the $\alpha_k$ are positive numbers summing to one. Since  $i' > i$
the point $x$ was not chosen in to be in $A_{i'j'}$ when it was defined. Part a) of the definition of $A_{i'j'}$
holds since $x \in A_{ij}$ and $F_{ij} \subset F = F_{i'j'}$, and part c) holds since $x$ made it to $A_{ij}$. As a result,
part b) of the definition fails; there is some $v_{k_0} \in v(f) \cap F_{i'j'}$ such that 
$$|x^{v_{k_0}}| < (C_{i'})^{-1} \sup_{v \in v(f) \cap F_{i'j'}} |x^v|$$
$$=   (C_{i'})^{-1}\sup_{v \in v(f) \cap F_{ij}}|x^v| \eqno (2.5)$$ 
The last equality follows from part a) of the definitions of $A_{ij}$ and $A_{i'j'}$. Note that $(2.4)$ can be rewritten as 
$$|x^w| = \prod_{v_k \in v(f) \cap F_{i'j'}}|x^{v_k}|^{\alpha_k} \eqno (2.6)$$
So in view of $(2.5)$ we have
$$|x^w| < (C_{i'})^{-{1 \over \alpha_{k_0}}}(\sup_{v \in v(f) \cap F_{ij}}|x^v|)^{\alpha_{k_0}}\prod_{v \in v(f), v \neq v_{k_0}}(\sup_{v \in v(f) \cap F_{ij}}|x^v|)^{\alpha_k} \eqno (2.7)  $$
Since the $\alpha_k$ sum to 1 the right-hand side of $(2.7)$ is
$$(C_{i'})^{-{1 \over \alpha_{k_0}}}\sup_{v \in v(f) \cap F_{ij}}|x^v|\eqno (2.8)$$
Since $x \in A_{ij}$, the above is at most
$$\leq (C_{i'})^{-{1 \over \alpha_{k_0}}}C_i\inf_{v \in v(f) \cap F_{ij}}|x^v|\eqno (2.9)$$
Since $i' > i$, if the constants $C_1,C_2,...,C_n$ were chosen to be increasing fast enough then for a uniform constant $\mu_0$ we have $(C_{i'})^{-{1 \over \alpha_{k_0}}}C_i \leq C_{i'}^{-\mu_0} \leq C_{i+1}^{-\mu_0}.$
Thus $(2.9)$ completes the proof of Lemma 2.0.\hfill$\square$

We now come to the chief result of this section, which gives a useful description of the way in which the terms $|x^v|$ for $v \in v(f) \cap F_{ij}$ dominate on $A_{ij}$. It is the generalization of Lemma 2.0 of [G3] to arbitrary $K$, with a streamlined proof.

\noindent {\bf Theorem 2.1.} If the number $N$ in the definition of the $C_i$ is chosen sufficiently large, the following hold
for all $i > 0$.

\noindent a) If $x \in E$ such that the following two conditions hold, then $x \in A_{ij}$.

\noindent 1) $\inf_{v\in v(f) \cap F_{ij}}|x^{v}| \geq C_i^{-1}\sup_{v \in v(f)}|x^v|$.  \hfill \break
\noindent 2) $\sup_{v \in v(f): v \notin F_{ij}}|x^{v}| < C_{n-1}^{-1} \inf_{v \in  v(f) \cap F_{ij} }|x^v|$.

\noindent b) There is a $\mu > 0$ depending only on $N(f)$ such that if $x \in A_{ij}$ then the following two conditions hold.

\noindent 1) $\inf_{v\in v(f) \cap F_{ij}}|x^{v}| \geq C_i^{-1}\sup_{v \in v(f)}|x^v|$.\hfill\break
\noindent 2) $\sup_{v \in v(f): v \notin F_{ij}}|x^{v}| \leq C_{i+1}^{-\mu} \inf_{v \in v(f) \cap F_{ij}}|x^v|$.

\noindent In the case where $i = 0$, we just have parts 2) of both a) and b). 

\noindent {\bf Proof.} Part a) is straightforward: If $x \in E$ satisfies the conditions of part a), then if $F_{i'j'}$ is such that
$i' > i$ or $i' = i$ and $j' < j$ then there is at least one vertex $v'$ such that $v' \in F_{i'j'}$ but $v' \notin F_{ij}$. As a result, 
for all $v \in F_{ij}$ one has $|x^{v'}|< C_{n-1}^{-1} |x^v| \leq C_{j'}^{-1}|x^v|$ . So condition b) of the definiton of $A_{i'j'}$ will
not be satisfied. Thus $x \notin A_{i'j'}$. Conditions 1) and 2) now ensure that $x$ will be selected to be a member of $A_{ij}$.
This completes the proof of part a).

\noindent We proceed to the more difficult part b). Note that condition 1) automatically holds by construction, so our goal is to 
show that condition 2) holds for all $x \in A_{ij}$. Let $F_{ij}$ be any 
compact face of $N(f)$, let $v \in v(f)$ be on $F_{ij}$ and let $v' \in v(f)$ not on $F_{ij}$. Our goal is to show that there
is a $\mu > 0$ independent of $(i,j)$ such that $|x^{v'}| < C_{i+1}^{-\mu}|x^v|$ whenever $x \in A_{ij}$.  

Let $F$ denote the
face of $N(f)$ of minimal dimension  containing $F_{ij}$ and $v'$; in the event there is no such face we define $F$ to be all of $N(f)$ which will serve as a substitute. 
Let $G$ denote the convex hull of $F_{ij}$ and $v'$. I claim that $G$  intersects the interior of $F$. To see this,  let
$\{w \in \R^n: a_j \cdot w = b_j \}$ be supporting hyperplanes of $N(f)$ transverse to $F$ such that the boundary of $F$ is given by
$\{w \in F: a_j \cdot w = b_j$ for some $j\}$; in the event that $F = N(f)$ we take these to be all supporting hyperplanes of $N(f)$.  Note that for each $j$, $\{w \in G: a_j \cdot w = b_j\}$ is either all of $G$ or a subset of $G$ of measure zero. It will not be
all of $G$, for if it were then $G$ would be a subset of $\{w \in F: a_j \cdot w = b_j\}$, a face of $N(f)$ of lower dimension than
$F$, and this would contradict the minimality of the dimension of $F$. We conclude that each $\{w \in G: a_j \cdot w = b_j\}$ is a subset
of $G$ of measure zero. Thus so is their union. In particular there exists some $w \in G$ with $a_j \cdot w \neq b_j$ for each $j$. This is
equivalent to $w$ being in the interior of $F$ as needed. 

From now on we view $w$ as fixed for a given $F_{ij}$ and $v'$; we will be applying Lemma 2.0 and we want the constants $\mu_0$ and
$N$ to depend on $N(f)$  only and not the choices of $w$.

Let $x \in A_{ij}$.  Without loss of
generality we may assume that $w \notin F_{ij} \cup \{v'\}$ since $F_{ij} \cup \{v'\}$ is a subset of the convex hull of $F_{ij}$
and $\{v'\}$ of measure zero. So $w$ is in the interior of the segment connecting 
$v'$ to some $a \in F_{ij}$ and we may write $w = \beta v' + (1 - \beta )a$ for some $0 < \beta < 1$. Thus we have
$$|x^{v'}| = |x^w|^{1 \over \beta} |x^{a}|^{-{1 - \beta \over \beta}}\eqno (2.10)$$
Using Lemma 2.0, this is at most
$$C_{i+1}^{-{\mu_0 \over \beta}}\big(\inf_{v \in v(f) \cap F_{ij}} |x^v|^{1 \over \beta}\big)
|x^{a}|^{-{1 - \beta \over \beta}} \eqno (2.11)$$
Since $a \in F_{ij}$, one may write $a$ as a convex combination $\sum_{v_k \in v(f) \cap F_{ij}}\gamma_k v_k$. Here the 
$\gamma_k$ are nonnegative numbers whose sum is 1. Thus we have
$$|x^{a}| = \prod_{v_k \in v(f) \cap F_{ij}}|x^{v_k}|^{\gamma_k}\eqno (2.12) $$
Because each $|x^{v_k}|$ in the above product is at least $\inf_{v \in v(f) \cap F_{ij}}|x^v|$, we have
$|x^{a}| \geq \inf_{v \in v(f) \cap F_{ij}}|x^v|$. Substituting this back into $(2.11)$ gives
$$|x^{v'}| \leq C_{i+1}^{-{\mu_0 \over \beta}}\big(\inf_{v \in v(f) \cap F_{ij}} |x^v|\big)^{1 \over \beta}
\big(\inf_{v \in v(f) \cap F_{ij}} |x^v|\big)^{-{1 - \beta \over \beta}}\eqno (2.13)$$
$$= C_{i+1}^{-{\mu_0 \over \beta}}\inf_{v \in v(f) \cap F_{ij}} |x^v| \eqno (2.14)$$
Since $v'$ was an arbitary vertex of $N(f)$ not on $F_{ij}$, setting $\mu$ to be the minimal ${\mu_0 \over \beta}$ we get part b) of Theorem 2.1 and we are done.\hfill$\square$

The next lemma is a way of describing how on $A_{ij}$, the terms $f_{\alpha}x^{\alpha}$ for $\alpha \in F_{ij}$ dominate the Taylor expansion of $f(x)$. For the case where
$K = \R$, it is Lemma 2.1 of [G3]. The proof there transfers over word for word to general $K$, so we will not reprove it here.

\noindent {\bf Lemma 2.2.} For any positive integer $d$, there is a constant $E_{d,f}$ depending on $d$ and $f(x)$
and a positive constant $\eta$ depending on $N(f)$ such that on $A_{ij}$ we have
$$\sum_{\alpha \notin F_{ij}} |f_{\alpha}||\alpha|^d |x^{\alpha}| < E_{d,f}(C_{i+1})^{-\eta}\sup_{v \in v(f)}|x^v|$$

\noindent {\heading Examples.}

\noindent {\bf Example 1.} 
Suppose $f(x,y)$ is a function of two real variables. Let $e$ be an edge of $N(f)$. Then $e$
joins two vertices of $N(f)$, which we denote by $v_1 = (\alpha_1,\beta_1)$ and $v_2 = (\alpha_2,\beta_2)$, where $\alpha_1 > \alpha_2$ and $\beta_1 
< \beta_2$. Let $-{1 \over m}$ be the slope of this edge. Then $\alpha_1  + m\beta_1 = \alpha_2 + m \beta_2$. So whenever $|y| 
= |x|^m$ one has that $|x^{\alpha_1}y^{\beta_1}| = |x^{\alpha_2}y^{\beta_2}|$. More generally, $|x^{\alpha_1}y^{\beta_1}|$ and
 $|x^{\alpha_2}y^{\beta_2}|$ are within a factor of $C_1$ of each other if $(x,y)$ is in the domain $B_e$ defined by
$$B_e = \{(x,y):  C_1^{-{1 \over \beta_2 - \beta_1}}|x|^m < |y| < C_1^{{1 \over \beta_2 - \beta_1}}|x|^m\}$$
Furthermore, for any other vertex $v = (\alpha,\beta)$ of $N(f)$, if $|y| = C|x|^m$
then $|x^{\alpha}y^{\beta}| = C^{\beta}|x|^c$ for some $c > \alpha_1 + m\beta_1 = \alpha_2 + m\beta_2$. Thus
if we are in a small enough neighborhood of the origin, then on $B_e$ each such $|x^{\alpha}y^{\beta}|$ will be less than both $C_1^{-1}|x^{\alpha_1}y^{\beta_1}|$ and $C_1^{-1}|x^{\alpha_2}y^{\beta_2}|$. Hence if $E$ denotes the square $\{(x,y) \in \R^2: 0 < |x|, |y| < C_2^{-1}\}$ for large enough $C_2$, we may take our $A_{1j}$ to be
the sets $B_e \cap E$ for the various edges $e$ of $N(f)$.

As for the $A_{0j}$,  between any two $A_{1j}$ corresponding to intersecting edges there will be an $A_{0j}$ of the form $\{(x,y) \in E: C_1^{1 \over {\beta_2 - \beta_1}}|x|^m < |y| < C_1^{-{1 \over \beta_2' - \beta_1'}}|x|^{m'}\}$. There will be two additional $A_{0j}$'s, one 
of the form
$\{(x,y) \in E: C_1^{{1 \over \beta_2 - \beta_1}}|x|^m < |y| \}$ and  the other of the form $\{(x,y) \in E:|y| < C_1^{-{1 \over \beta_2 - \beta_1}}|x|^{m}\}$, corresponding the uppermost and lowermost vertices of $N(f)$ respectively. We are excluding
the case where $N(f)$ has only one vertex, in which case there are no $A_{1j}$'s and exactly one $A_{0j}$, consisting of all of $E$.

\noindent {\bf Example 2.} 
Suppose now we are in three dimensions, again in the real case. Let $F$ be a two-dimensional face of $N(f)$.
 If $n = (n_1,n_2,n_3)$ denotes a vector normal to this face with each $n_i > 0$, then for $v \in v(f)$, $n \cdot v$ is maximized for $v \in F$.
 In particular, on a curve $\{(c_1 t^{n_1}, c_2t^{n_2},c_3t^{n_3}): t > 0\}$, the functions $|x^v|$ for $v \in F$ dominate
the $|x^v|$ for $v \notin F$ as $t \rightarrow 0$, while the $x^v$ for $v \in F$ remain within a fixed factor of each other. 
Thus by Theorem 2.1 the $A_{2j}$ corresponding to $F$ will be a "horn" approximately consisting of a union of 
such curves with positive $c_1,c_2,$ and $c_3$ within some power of $C_2$ of each other, along with its 7 reflected horns about the 
coordinate planes.

Now let $e$ be a one-dimensional edge of $N(f)$, joined by two vertices $v_1$ and $v_2$. By Theorem 2.1, the corresponding $A_{1j}$ will contain  points in a small cube centered at the origin for which $|x^{v_1}|$ and $|x^{v_2}|$ are within a factor of $C_1$ of each other,
such that $|x^v|$ for vertices $v$ not on $e$ are of magnitude  less than $C_2^{-\mu}$ times the smaller of $|x^{v_1}|$ or $|x^{v_2}|$.
It will contain all points in this cube for which $|x^{v_1}|$ and $|x^{v_2}|$ are within $C_1$ of each other,
and  for which $|x^v|$ for vertices $v$ not on $e$ are of magnitude  less than $C_2^{-1}$ times the smaller of $|x^{v_1}|$ or $|x^{v_2}|$. In the event that there are no vertices other than $v_1$ and $v_2$, then $A_{1j}$ consists of all points in the cube for 
which $|x^{v_1}|$ and $|x^{v_2}|$ are within $C_1$ of each other, similar to the situation for the $A_{1j}$ in the
2-dimensional case. It can be shown that for the vertices not on $e$ one may restrict consideration to vertices $v$ not on $e$ which are on two-dimensional compact faces of $N(f)$ containing $e$. 

 Geometrically, the $A_{1j}$ again can be viewed as horns in the upper right quadrant along with 7 reflections 
about the coordinate planes, but this
time the horns will be not be contained in any bundle of curves in the sense the $A_{2j}$ were.

Lastly if $v$ is a vertex of $N(f)$, then by Theorem 2.1 the corresponding
$A_{0j}$ will consist of points in a small cube centered at the origin for which $|x^v| > C_1^{\mu}|x^{v'}|$ for all other vertices $v'$ of $N(f)$, and will contain all points in this cube for which $|x^v| > C_2|x^{v'}|$ for all vertices $v'$ of $N(f)$. Similar to 
above it can be shown that one can restrict to $v'$ on edges of $N(f)$ containing $v$. In the event that $v$ is the only vertex
of $N(f)$, then $A_{01}$ will be an entire cube centered at the origin, and there will be no other $A_{ij}$.

\noindent {\bf Example 3.} Consider the case in $n$ dimensions where $N(f)$ has a vertex on each coordinate axis at some height $h$, the same for
each axis, and no other vertices. Then the $A_{ij}$ are as follows. For each nonempty $I \subset \{1,...,n\}$ there is a face $F_{ij}$ for $i = |I|$
which is the convex hull of the vertices on the $x_i$ axis for $i \in I$. By Theorem 2.1, the corresponding $A_{ij}$ consists of points
in a small cube centered at the
origin for which the $|x_i|$ for $i \in I$ are within a factor of $C_i^{1 \over h}$ of the maximal $|x_i|$, achieved for some $i \in I$,  and such that $|x_i|$ for $i \notin I$ are 
less than $C_{i+1}^{-{\mu \over h}}|x_i|$ for all $i \in I$. Furthermore, $A_{ij}$ contains all the points in this cube for which $|x_i|$ for $i \in I$ are within a factor of $C_i^{1 \over h}$ of the maximal $|x_i|$, and such that $|x_i|$ for $i \notin I$ are 
less than $C_{n-1}^{-{\mu \over h}}|x_i|$ for all $i \in I$.

\noindent {\heading 3. Performing the blowups; the $A_{ijk}'$ in the new coordinates.} 

Our first lemma gives a special case of resolution of singularities which will be used as a building block for the resolution of 
singularities algorithm of this paper. Let $K_0 = K - \{0\}$ as before. In our first lemma, one starts with the cube $F =\{ x \in K_0^n: |x_l| \leq 1$ for all $l\}$, divides it
into sets  $F_1 = \{x \in F: |x_j| \leq |x_k|\}$  and $F_2 = \{x \in F: |x_j| > |x_k|\}$ for some variables $x_j$ and $x_k$. One then
blows up the sets $F_1$ and $F_2$, changing variables $x_j = x_j x_k$ in $F_1$ and $x_k = x_kx_j$ in $F_2$, leaving the other
variables unchanged. The result are two cubes $F_1'$ and $F_2'$ (not necessarily closed this time), which we then 
blow up in an analogous fashion. This will be repeated as needed to obtain the lemma.

\noindent {\bf Lemma 3.1.} Let $F =\{ x \in K_0^n: |x_l| \leq 1$ for all $l\}$, and
let $m_1(x),...,m_p(x)$ be distinct monomials. Then there is a way of starting with $F$ and then doing a finite sequence of blowups in the above sense  after which one has  cubes $E_1,...,E_q$ (now not necessarily closed) such that in the blown-up coordinates of a given cube $E_i$, the monomials $m_1(x),...,m_p(x)$ are new monomials $\tilde{m}_1(x),...,\tilde{m}_p(x)$ whose exponents are linearly ordered. That is, if $\tilde{m}_i(x) = x^{\alpha}$
and $\tilde{m}_j(x) = x^{\beta}$, then either $\beta_k \geq \alpha_k$ for each $k$ or $\beta_k \leq \alpha_k$ for each $k$.

\noindent {\bf Proof.}  We can write the collection of all possible ${m_i(x) \over m_j(x)}$ as
$\{{x^{\alpha_{ij}} \over x^{\beta_{ij}}}: i,j = 1,...,p\}$ where for a given $i$ and $j$ the monomials $x^{\alpha_{ij}}$ and $x^{\beta_{ij}}$ have no 
variable in common. In order to prove the lemma, we need to perform a sequence of blowups after which on a given $E_k$, the blown up  ${x^{\alpha_{ij}} \over x^{\beta_{ij}}}$ becomes either a monomial or a reciprocal of a monomial. If one 
accomplishes this for one  ${x^{\alpha_{ij}} \over x^{\beta_{ij}}}$, this function will still be a  monomial or a reciprocal of a monomial 
respectively after any future blowups. So it suffices to prove we can do this for a single ${x^{\alpha_{ij}} \over x^{\beta_{ij}}}$ as
we can then just repeat the process for each  ${x^{\alpha_{ij}} \over x^{\beta_{ij}}}$ one after another.

So we suppress the indices $i$ and $j$ and focus our attention on a single function ${x^{\alpha} \over x^{\beta}}$, for which we will
show there exists the desired sequence of blowups. Writing $\alpha = (\alpha_1,...\alpha_n)$ and $\beta = (\beta_1,...,\beta_n)$, we let $a = \max_m \alpha_m$, and let $b = \max_m \beta_m$. The proof is first by induction on $c =\max(a,b)$. Given $c$ we induct on the number $M$ of variables $x_i$ appearing to the $c$ power in $x^{\alpha}$ or $x^{\beta}$ (To be precise, in our arguments either $M > 1$ and we will decrease $M$ by $1$, or $M = 1$ and we decrease the value of $c$.) The induction starts with
the trivial $(c,M) = (1,1)$ case.

\noindent {\bf Case 1. a = b}. Suppose $a = b = c$. Let $x_j$ be a variable appearing in $x^{\alpha}$ to the $c$th power,  and let
$x_k$ be a variable appearing in $x^{\beta}$ to the $c$th power. If one performs a blowup in the $x_j$ and $x_k$ variables, one
obtains two regions. On the first region $x_j$ becomes $x_jx_k$ with the other variables remaining fixed, and on the second region $x_k$ becomes $x_jx_k$ with the other variables remaining fixed. On the first region the $x_j^c$ in the numerator becomes 
$x_j^cx_k^c$, and the new $x_k^c$ appearing in the numerator can be cancelled with the $x_k^c$ appearing in the denominator.
Similarly, in the second region, an extra $x_j^c$ appears in the denominator which can be cancelled with the $x_j^c$ in the
numerator. Thus in either case, one of the variables appearing to the $c$th power in the overall ratio ${x^{\alpha} \over x^{\beta}}$ dissappears and the expression otherwise is unchanged. Thus we have decreased $M$ by 1 and we may now apply
the induction hypothesis. Thus we are done with the argument for when $a = b$.

\noindent {\bf Case 2. a $\neq$ b}. Our proof is symmetric in $\alpha$ and $\beta$, so without loss of generality we assume that $a > b$. Given our $(c,M) (= (a,M)$), we induct on the number $N$ of indices $i$ for which $x_i$ appears in $x^{\beta}$ to any
power. This final induction starts with the $N = 0$ case where $x^{\beta} = 1$ and the result is immediate.

Let $(a,M,N)$ be the case we are currently trying to prove. Let $x_j$ be a variable appearing to the $a$th power in $x^{\alpha}$, and let $x_k$ be any variable appearing in $x^{\beta}$. Define $d$ to be the number
such that $x_k$ appears to the $d$ power in $x^{\beta}$, where $d<a$. We perform a blowup in the $x_j$ and $x_k$ 
variables. We obtain two domains. On the first domain, the $x_j^a$ in the numerator becomes $x_j^ax_k^a$. The $x_k^a$ in
the numerator can be cancelled with the $x_k^d$ the denominator, resulting in an $x_k^{a-d}$ in the numerator and there is no longer the variable $x_k$ appearing in the denominator. Thus $N$ is decreased and by induction hypothesis we are done for the
first domain. On the second domain, the $x_k^d$ in the denominator becomes $x_j^dx_k^d$, and the $x_j^a$ in the
numerator can be cancelled with the $x_j^d$ in the denominator, obtaining a $x_j^{a-d}$ in the numerator in place of $x_j^a$, 
with the function  ${x^{\alpha} \over x^{\beta}}$ otherwise unchanged. Thus if $M > 1$, $M$ becomes $M-1$, and if $M = 1$ 
$c = a$ gets decreased. Thus by induction hypothesis in the second domain we again are done. This concludes the proof. 
\hfill$\square$

\noindent {\bf Blowing up the domains $A_{ij}$}

We now let $f(x) = \sum_{\alpha} c_{\alpha} x^{\alpha}$ be a the power series of an analytic function on a neighborhood of the origin. We apply Lemma 3.1 to the monomials $\{x^v: v$ is a vertex of the Newton polygon $N(f)\}$, obtaining the resulting
collection of cubes, which we denote by
$E_1,...,E_q$. Let $\gamma_k$ denote the composition of the blowups used in creating $E_k$. Then there is 
a linear function $L_k$ such that the coordinate change $\gamma_k$ converts these $x^v$
 into the monomials $x^{L_k(v)}$  and furthermore the exponents $\{L_k(v): v $ is a vertex of $N(f)\}$ are linearly ordered under the lexicographic ordering by Lemma 3.1.

We next examine what happens to the domains $A_{ij}$ from Section 2 under these sequences of blowups. Define $A_{ijk} 
= A_{ij} \cap \gamma_k(E_k)$. Since the $A_{ij}$ were (up to a set of measure zero) a subdivision of  the cube $E = \{x \in  K_0^n: |x_l| < C_n^{-1}$ for all $l \}$, the $A_{ijk}$ are also a subdivision of $E$. Define $A_{ijk}' =\gamma_k^{-1}(A_{ijk}) = \gamma_k^{-1}(A_{ij}) \cap E_k$, and define $G_k = \gamma_k^{-1}(E) \cap E_k$. Then for a fixed $k$, the $A_{ijk}'$ are a subdivision of $G_k$. Note that if we write $\gamma_k(x) = (\gamma_{k1}(x),...,
\gamma_{kn}(x))$, then we have
$$G_k= \{x \in  E_k: |\gamma_{kl}(x)| < C_n^{-1} {\rm\,\,for\,\,all\,\,} l\}$$
Next, we examine how Theorem 2.1 translates in the blown up coordinates for a given $k$. $G_k \subset E_k \subset\{x \in K_0^n: |x_l| \leq 1$
for all $l\}$, so since the exponents of the transformed monomials $x^{L_k(v)}$ are linearly ordered, there is a corresponding ordering of the functions $x^{L_k(v_i)}$; for a given $p$ and $q$ either $|x^{L_k(v_p)}| \leq |x^{L_k(v_q)}|$ throughout
$G_k$, or $|x^{L_k(v_q)}| \leq |x^{L_k(v_p)}|$ throughout $G_k$.

This makes translating Theorem 2.1 in the blown up coordinates rather straightforward. As in section 2, we let $v(f)$ denote the set of vertices of $N(f)$. Since the functions $\{|x^{L_k(v)}| : v \in v(f)\}$ are linearly ordered on $G_k$ in the sense of the above paragraph, there are vertices 
$v',v'',v'''$ of $N(f)$ such that throughout $G_k$ one has the following
$$|x^{L_k(v')}| = \sup_{v \in v(f)}|x^{L_k(v)}|$$
$$|x^{L_k(v'')}| = \inf_{v \in v(f) \cap F_{ij}}|x^{L_k(v)}|$$
$$|x^{L_k(v''')}| = \sup_{v \in v(f) \cap F_{ij}^c}|x^{L_k(v)}|$$
Note that $L_k(v') \leq L_k(v'')$ under the lexicographical ordering, with strict inequality if and only if $i > 0$. Note also that by Theorem 2.1 b) part 2, if $(i,j,k)$ is such that $A_{ijk}'$ is nonempty then we have $L_k(v'') < L_k(v''')$. Thus for $i > 0$ we may define the  nonconstant monomial  $p_{ijk}(x) = x^{L_k(v'') - L_k(v')}$ and for all $i$ we may define the nonconstant monomial $q_{ijk}(x)= x^{L_k(v''') - L_k(v'')}$. 
 Then Theorem 2.1 translates into the following.

\noindent {\bf Lemma 3.2.} If the number $N$ in the definition of the $C_i$ is chosen sufficiently large, then for all
$i > 0$ and all $j$ and $k$ we have (up to a set of measure zero)
$$\{x \in G_k: |p_{ijk}(x)| \geq C_i^{-1}{\rm\,\,and\,\,} |q_{ijk}(x)| < C_{n-1}^{-1}\} \subset A_{ijk}'$$
$$A_{ijk}' \subset \{x \in G_k:
|p_{ijk}(x)| \geq C_i^{-1} {\rm\,\,and\,\,} |q_{ijk}(x)| < C_{i+1}^{-\mu}\} \eqno (3.1a)$$
If $i = 0$, then for all $j$ and $k$ we have  (up to a set of measure zero)
$$\{x \in G_k: |q_{0jk}(x)| < C_{n-1}^{-1}\} \subset A_{0jk}' \subset \{x \in G_k: |q_{0jk}(x)| < C_{1}^{-\mu}\} \eqno (3.1b)$$
Since each $G_k$ is a subset of the unit cube, if $x \in G_k$ such that $|p_{ijk}(x)| \geq C_i^{-1}$, then $|x_r| \geq C_i^{-{1 \over
deg(p_{ijk})}}$ for all $r$ such that $x_r$ appears in $p_{ijk}(x)$. This fact will be highly relevant for our arguments. For example,
it implies that 
we can assume that not all variables appear in $p_{ijk}(x)$. For if they did, then on $A_{ijk}'$ we would have that $|x_r| > C_i^{-{1 \over deg(p_{ijk})}}$ for all $r$. Then assuming $C_{i+1}$ was chosen large enough relative to $C_i$, the right-hand side of $(3.1a)$
implies that $A_{ijk}'$ is empty.

Thus when $i > 0$,  permuting the variables if necessary we may let $a$ be such that the variables 
$x_1,...,x_p$ do not appear in $p_{ijk}(x)$ and the
variables $x_{p+1},...,x_n$ do appear in $p_{ijk}(x)$.  We rename our variables now, letting $y_l = x_l$ for $1 \leq l \leq p$, and
$z_l = x_{l + p}$ for $1 \leq l \leq n - p$. We correspondingly write $q_{ijk}(x)$ as $s_{ijk}(y)t_{ijk}(z)$ and $p_{ijk}(x)$ as
$p_{ijk}(z)$. For consistency of notation,  when $i = 0$ we rename each $x_l$ variable $y_l$, and the rename the function
$q_{0jk}(x)$ as $s_{0jk}(y)$. Then Lemma 3.2 can be recast as follows.

\noindent {\bf Lemma 3.2'}. If $i > 0$, then we have
$$\{(y,z) \in G_k: |s_{ijk}(y)| < C_{n-1}^{-1} {\rm\,\,and\,\,}  |p_{ijk}(z)| \geq C_i^{-1}\} \subset A_{ijk}'$$
$$A_{ijk}' \subset \{(y,z) \in G_k: |s_{ijk}(y)|< C_{i+1}^{-\mu}
{\rm\,\,and\,\,}|p_{ijk}(z)| \geq C_i^{-1}  \} \eqno (3.2a)$$
If $i = 0$ then we have
$$\{y \in G_k: |s_{0jk}(y)| < C_{n-1}^{-1}\} \subset A_{0jk}' \subset \{y \in G_k: |s_{0jk}(y)| < C_{1}^{-\mu}\} \eqno (3.2b)$$
\noindent {\bf Proof.} The $i = 0$ case is identical to that of Lemma 3.2, so we assume $i > 0$. $(3.1a)$ can be rewritten as
$$\{(y,z) \in G_k: |s_{ijk}(y)t_{ijk}(z)| < C_{n-1}^{-1} {\rm\,\,and\,\,}   |p_{ijk}(z)| \geq C_i^{-1}\} \subset A_{ijk}'$$
$$A_{ijk}' \subset \{(y,z) \in G_k:|s_{ijk}(y)t_{ijk}(z)|< C_{i+1}^{-\mu}
{\rm\,\,and\,\,}|p_{ijk}(z)| \geq C_i^{-1}  \} \eqno (3.3)$$
Since $C_i^{-{1 \over deg(p_{ijk})}} \leq |z_l| \leq 1$ for all $l$ on $A_{ijk}'$, on $A_{ijk}'$ one has $C_i^{-{deg(t) \over deg(p)}} \leq |t_{ijk}(z)| \leq 1$. Thus $(3.3)$ implies
$$\{(y,z) \in G_k: |s_{ijk}(y)| < C_{n-1}^{-1}{\rm\,\,and\,\,} |p_{ijk}(z)| \geq C_i^{-1} \} \subset A_{ijk}'$$
$$A_{ijk}' \subset \{(y,z) \in G_k: |s_{ijk}(y)|< C_i^{{deg(t) \over deg(p)}}C_{i+1}^{-\mu}
{\rm\,\,and\,\,}  |p_{ijk}(z)| \geq C_i^{-1} \} \eqno (3.4)$$
If the constants $C_i$ were chosen to increase rapidly enough, one has $C_i^{{deg(t) \over deg(p)}}C_{i+1}^{-\mu}
< C_{i+1}^{-{\mu \over 2}}$.  Thus replacing $\mu$ by ${\mu \over 2}$ if necessary equation $(3.2a)$ follows and we are done.
\hfill $\square$

Note that the right-hand sides of $(3.2a)-(3.2b)$ ensure that there is at least one $y_l$ variable appearing in $s_{ijk}(y)$ (i.e. no $s_{ijk}(y)$ is ever just $1$).

\noindent {\bf Lemma 3.3.} For each $(i,j,k)$ with $i > 0$ there is an $\alpha_{ijk}$ such that each $x^{L_k(v)}$ for a vertex $v \in F_{ij}$ can be written as $y^{\alpha_{ijk}}z^{\eta}$ for some $\eta$, and such that each $x^{L_k(v)}$ for a vertex $v \notin F_{ij}$ can be
written as $y^{\alpha_{ijk} '}z^{\eta}$ with $\alpha_{ijk}' > \alpha_{ijk}$ under the lexicographic ordering.

\noindent {\bf Proof.}
Since the $z$ variables were defined to be the variables for which the minimal $|x^{L_k(v)}|$ for a vertex $v \in F_{ij}$ 
differs from the maximal $|x^{L_k(v)}|$ over all $v$ (and therefore over all $v \in F_{ij}$ by part a) of the definition of the
 $A_{ij}$),  all $x^{L_k(v)}$ for $v \in F_{ij}$ will be the same in the $y$ variables. In other words, there's a
fixed $\alpha_{ijk}$ such that each $x^{L_k(v)}$ for $v \in F_{ij}$ can be written as $y^{\alpha_{ijk}}z^{L'(v)}$ for some linear
map $L'$. Furthermore, since $s_{ijk}(y)t_{ijk}(z)$ was defined as the ratio of the minimal $|x^{L_k(v)}|$ for $v \in F_{ij}$ (now known to be  of 
the form $y^{\alpha_{ijk}}z^{\eta}$ for some $\eta$) to the maximal 
$|x^{L_k(v)}|$ for $v \notin F_{ij}$ and $s_{ijk}(y)$ has at least one $y_l$ variable appearing by above, each $x^{L_k(v)}$ for
 $v \notin F_{ij}$ will also contain this $y_l$ variable and therefore can be 
written as $y^{\alpha_{ijk} '}z^{L'(v)}$ with $\alpha_{ijk}' > \alpha_{ijk}$ under the lexicographic ordering. This completes the proof. \hfill$\square$

\noindent {\heading Example. The two-dimensional case.}

\noindent Suppose $f(x,y)$ is a function of two variables. Then by Lemma 3.2', a given $A_{1jk}'$ satisfies
$$\{(x,y) \in G_k: |s_{1jk}(x)| < C_1^{-1}{\rm\,\,and\,\,} |p_{1jk}(y)| \geq C_1^{-1}\} \subset A_{1jk}'$$
$$A_{1jk}' \subset \{(x,y) \in G_k: |s_{1jk}(x)|< C_2^{-\mu}
{\rm\,\,and\,\,}  |p_{1jk}(y)| \geq C_1^{-1} \} \eqno (3.5)$$
The structure of the $G_k$ provides restrictions on what an $A_{1jk}'$ satisfying $(3.5)$ can be. Since we are in two dimensions, up to a set of measure zero $G_k$ is given by
$$G_k = \{(x,y): |x|, |y| < 1,\,\,|\gamma_{k1}(x,y)| < C_2^{-1},\,\, |\gamma_{k2}(x,y)| < C_2^{-1}\} \eqno (3.6)$$
Here $\gamma_{k1}(x,y)$ and $\gamma_{k2}(x,y)$ are monomials and are
the components of the composition of blowups $\gamma_k(x,y)$. If either $\gamma_{k1}(x,y)$ or $\gamma_{k2}(x,y)$ is a function
of $y$ only, then $G_k \subset \{(x,y): |y| > C_2^{-{1 \over m}}\}$ for some $m$, and if the $C_i$ were
chosen to be increasing fast enough any set $A_{1jk}'$ satisfying $(3.5)$ would be empty. So $\gamma_{k1}$ and $\gamma_{k2}$ 
must depend on $x$ as well as possibly $y$. 

Since $p_{1jk}(y)$ is a monomial, if $y$ is such that $|p_{1jk}(y)| > C_1^{-1}$, then $|y| > C_1^{-{1 \over deg(p)}}$. Then
a condition $|x^ay^b| < C_2^{-1}$ with  $a> 0$, $b \geq 0$ implies that $|x| < C_2^{-{1 \over a}}C_1^{b \over a deg(p)}$. Assuming the
$C_i$ were chosen to increase fast enough this means that on $G_k \cap A_{1jk}'$ for some $a$ we have $|x| < C_2^{-{1 \over 2a}}$. In this
case the condition in $(3.5)$ that  $|s_{1jk}(x)| < C_1^{-1}$ or $C_2^{-\mu}$ are superfluous, and so whenever $i = 1$ we have
$$A_{1jk}' = \{(x,y) \in G_k: |p_{1jk}(y)| \geq C_1^{-1}\} \eqno (3.7)$$
Since $p_{1jk}(y)$ is just a monomial, this is the same as
$$A_{1jk}' = \{(x,y) \in G_k: |y| \geq C_1^{-{1 \over deg(p)}}\} \eqno (3.7')$$
Note that by above, $|x| < C_2^{-{1 \over 2a}}$ on $A_{1jk}'$. 

\noindent We now look at the sets $A_{0jk}'$. By Lemma 3.2' we have
$$\{(x,y) \in G_k: |s_{0jk}(x,y)| < C_1^{-1}\} \subset A_{0jk}' \subset \{(x,y) \in G_k: |s_{0jk}(x,y)| < C_1^{-\mu}\} \eqno (3.8)$$
If $s_{0jk}(x,y)$ contains both the $x$ and $y$ variables, then by $(3.6)$ and the left-hand side of $(3.8)$ the set $A_{0jk}'$ will be all of $G_k$. If 
$s_{0jk}(x,y)$ is a function of just one variable, say $y$, there are a few possibilities. If one or both of $\gamma_{k1}$ and $\gamma_{k2}$ is a 
function of $y$ only, then by $(3.6)$ $A_{0jk}'$ will be all of $G_k$ like above. If one or both are functions of $x$ only, then $G_k \subset
\{(x,y): |x| < C_2^{-d}, |y| < 1\}$ for some $d > 0$, so by $(3.8)$ $A_{0jk}'$ will be contained in the small box $\{(x,y): |x| < C_2^{-d}, |y| < C_1^{-{\mu \over deg(s_{0jk})}}\}$. If  $\gamma_{k1}$ and $\gamma_{k2}$ are both functions of $x$ and $y$, then $A_{0jk}'$ is a 
short but wide region of $|y|$ height at most $C_1^{-{\mu \over deg(s_{0jk})}}$ but which extends in the $x$-direction all the way to $|x| = 1$.

 Note that due to the forms we have determined for the various $A_{ijk}'$, if $G_k - A_{0jk}'$ is nonempty it must consist of a single
$A_{1j'k'}'$. 

In preparation for our next lemma, for a given face $F_{ij}$ of $N(f)$ write $f(x) = f_{F_{ij}}(x) + g_{ij}(x)$  (cf Definition 2.2). Recalling that $\gamma_k$ was defined to be the composition of blowups corresponding to $E_k$, we let $u_{ijk} = f_{F_{ij}} \circ \gamma_{k}$  and let $v_{ijk} =  g_{ij} \circ \gamma_{k}$. It is a consequence of Lemma 3.3 that when $i > 0$ each monomial of the polynomial $u_{ijk}(y,z)$ is of the form $y^{\alpha_{ijk}}z^{\eta}$ and each monomial of the Taylor series of $v_{ijk}(y,z)$ at the origin is of
 the form $y^{\alpha_{ijk}'}z^{\eta}$ for some $\alpha_{ijk}' \geq  \alpha_{ijk}$. When $i = 0$ the same holds if we interpret
$z^{\eta}$ to be $1$. In particular, when $i > 0$ we can write $u_{ijk}(y,z) = y^{\alpha_{ijk}}
 U_{ijk}(z)$ and  $v_{ijk}(z) = y^{\alpha_{ijk}} V_{ijk}(y,z)$ for analytic functions
$U_{ijk}$ and $V_{ijk}$, so that $f \circ \gamma_{k}(y,z) = y^{\alpha_{ijk}}F_{ijk}(y,z)$ where $F_{ijk}(y,z) = U_{ijk}(z) +  V_{ijk}(y,z)$. When $i = 0$ the 
analogous statement is that  $f \circ \gamma_{0jk}(y)$ can be written in the form $y^{\alpha_{0jk}}F_{0jk}(y)$ for some
analytic $F_{0jk}(y)$ with nonzero constant term.

Our next lemma shows that if $o_{ij}$ is
the maximum order of any zero of $f_{F{ij}}(x)$  on $K_0^n$, then if $i > 0$ on $A_{ijk}'$ we have corresponding lower bounds on 
some derivative of $f \circ \gamma_k(y,z)$  in the $z$ variables of order at most $o_{ij}$, with an analogous statement in the 
$i = 0$ case.

\noindent {\bf Lemma 3.4.} Assuming the constants $C_i$ were chosen to be increasing
sufficiently quickly the following hold.

\noindent {\bf a)} Suppose $i > 0$.  Let $o_{ij}$ be the maximum order of any zero of $f_{F_{ij}}(x)$ on 
$K_0^n$. Let $\pi(y,z) = z$ denote projection onto the $z$ coordinates, and let $B_{ijk} = \pi(A_{ijk}')$. Then $B_{ijk}$ 
can be written as a finite union $B_{ijk} =  \cup_lB_{ijkl}$ such that for each $l$ there is a 
directional derivative $\sum_{m=1}^{n-p}\beta_m \partial_{z_m}$ with each $\beta_m$ rational
and $\sum_m |\beta_m| = 1$, a $0 \leq q \leq o_{ij}$, and a constant $\delta_{ijkl} > 0$ such that for all $(y,z) \in A_{ijk}'$ with $z \in B_{ijkl}$ we have
$$  |(\sum_m \beta_m \partial_{z_m})^q F_{ijk}(y,z)| 
\geq \delta_{ijkl}\eqno (3.9a)$$
 Equivalently, $ |(\sum_m \beta_m \partial_{z_m})^q (f\circ \gamma_k)(y,z)|
\geq \delta_{ijkl} |y^{\alpha_{ijk}}| $.
Here $\alpha_{ijk}$ is as in Lemma 3.3.

\noindent {\bf b)} Suppose $i = 0$, so that $F_{0j}$ is a vertex $v$ of $N(f)$. Let $\alpha_{0jk}$ be such that $x^v$ becomes $y^{\alpha_{0jk}}$ in the blown-up coordinates. Then there is a constant $\delta_{0jk} > 0$ such that for all  $y \in A_{0jk}'$ we have
$$   | F_{0jk}(y)| \geq \delta_{0jk}\eqno (3.9b)$$
Equivalently, $| f \circ \gamma_{k}(y)| \geq \delta_{0jk} |y^{\alpha_{0jk}}| $.

\noindent {\bf Proof.} 
When $i > 0$, the condition that all
zeroes of $ f_{F_{ij}}(x) $ on $K_0^n$ have order at most $o_{ij}$ implies that all zeroes of $f_{F_{ij}} \circ \gamma_k(y,z) = u_{ijk}(y,z) = y^{\alpha_{ijk}} U_{ijk}(z)$ on $K_0^n$ also have order at most $o_{ij}$. Thus the same is true for $U_{ijk}(z)$.
As is well known (see Ch 8 of [St] for details), the space of partial derivative operators $\partial^{\alpha}$ of order $q$ is spanned by a finite list of directional derivative operators $(\sum_m \beta_m \partial_{z_m})^q$, where each
$\beta_m$ is rational. Without loss of generality, we may assume $\sum_m |\beta_m| = 1$. Thus when $i > 0$ one can write 
$B_{ijk}$ as the finite union $ \cup_lB_{ijkl}$ such that for each $l$ there is a 
directional derivative $\sum_m \beta_m \partial_{z_m}$, a $0 \leq q \leq o_{ij}$, and a
$\delta_{ijkl} > 0$ such that for all $z \in B_{ijkl}$ we have
$$|(\sum_m \beta_m \partial_{z_m})^q U_{ijk}(z)| > 2\delta_{ijkl}  \eqno (3.10)$$
As a result, for $(y,z)$ with $z \in B_{ijkl}$ we have
$$|(\sum_m \beta_m \partial_{z_m})^q u_{ijk}(y,z)| \geq 2\delta_{ijkl}
 |y^{\alpha_{ijk}}| \eqno (3.11)$$
If $i = 0$, $(3.11)$ still holds if we take $q = 0$ since $u_{0jk}(y)$ is just a monomial $U_{0jk}y^{\alpha_{0jk}}$ here. 

Write the Taylor expansion of $v_{ijk}(y,z)$ as $\sum_{\gamma,\delta} g_{\gamma \delta}y^{\gamma}z^{\delta}$ (Only the
$y$ variable appears if $i = 0$). When applying the derivative 
 $(\sum_m \beta_m \partial_{z_m})^q$ to a given monomial $ g_{\gamma \delta}y^{\gamma}z^{\delta}$,
one obtains at most $n^q$ terms, each having a coefficient with magnitude at most $(|\gamma| + |\delta|)^q|g_{\gamma \delta}|$.
 Since each of the $|z_m|$ are bounded below by $C_i^{-1}$ on $A_{ijk}'$, the magnitude of the 
corresponding term is at most $C_i^q(|\gamma| + |\delta|)^q|g_{\gamma \delta}y^{\gamma}z^{\delta}|$. Thus the
magnitude of the sum of all of these at most $n^q$ terms is at most $C_i^q n^q(|\gamma| + |\delta|)^q|g_{\gamma \delta}y^{\gamma}z^{\delta}|  \leq C_i^{o_{ij}}n^{o_{ij}}(|\gamma| + |\delta|)^{o_{ij}}|g_{\gamma \delta}y^{\gamma}z^{\delta}| $. 

Translating this back into the original $x$-coordinates, if $f_{\alpha}x^{\alpha}$ denotes the term of $f(x)$'s
Taylor expansion corresponding to $g_{\gamma \delta} y^{\gamma}z^{\delta}$, on $A_{ijk}'$ we have
$$|(\sum_m \beta_m \partial_{z_m})^q g_{\gamma \delta} y^{\gamma} z^{\delta}| \leq C_i^{o_{ij}}n^{o_{ij}}(|\gamma| + |\delta|)^{o_{ij}}|f_{\alpha}x^{\alpha}|$$
 Because the transition from $x$ to $(y,z)$ coordinates is via
the monomial map $\gamma_k$, there exists a constant $D_k$ such that $|\gamma| + |\delta|\leq D_k|\alpha|$. Therefore on $A_{ijk}'$ we have
$$|(\sum_m \beta_m \partial_{z_m})^q g_{\gamma \delta} y^{\gamma}z^{\delta}| \leq  (C_i D_k n)^{o_{ij}}|\alpha|^{o_{ij}}|f_{\alpha} x^{\alpha}| \eqno (3.12)$$
Since the $g_{\gamma \delta} y^{\gamma}z^{\delta}$ are the terms of $v_{ijk}(y,z)$'s Taylor expansion, by definition of 
$v_{ijk}(y,z)$ the terms $f(x)$'s Taylor expansion corresponding to the $g_{\gamma \delta} y^{\gamma}z^{\delta}$ are
exactly those with $\alpha \notin F_{ij}$. Thus adding $(3.12)$ over all $(\gamma,\delta)$ we get
$$|(\sum_m \beta_m \partial_{z_m})^qv_{ijk}(y,z)| \leq  (C_i D_{k}n)^{o_{ij}}\sum_{\alpha \notin F_{ij}} 
|\alpha|^{o_{ij}}|f_{\alpha}||x^{\alpha}| \eqno (3.13)$$
By Lemma 2.2, for some $\eta > 0$ the right-hand side of $(3.13)$ is bounded by $C_i^{o_{ij}}$ times a constant depending on $f$ times $C_{i+1}^{-\eta}
\sup_{v \in v(f)}|x^v|$. Note that since $x \in A_{ij}$, this supremal $|x^v|$ will occur for some $v \in F_{ij}$, so  $|x^v|$ can be written as  $|y^{\alpha_{ijk}}z^{\xi}|$ for some $\xi$. This in turn is at most $|y^{\alpha_{ijk}}|$ since
each $|z_m| \leq 1$. So for a constant $M$ depending on $f$ we have
$$|(\sum_m \beta_m \partial_{z_m})^q v_{ijk}(y,z)| \leq M C_i^{o_{ij}}C_{i+1}^{-\eta}|y^{\alpha_{ijk}}|$$
By $(3.10)$, $\delta_{ijkl}$ is bounded below by the minimum absolute value of certain derivatives of $U_{ijk}(z)$ on $\{z \in K^{n-p}: C_i^{-1} \leq |z_l| \leq 1$ for all $l\}$. So as long as
the constants  $C_i$ were chosen to be increasingly sufficiently quickly, for $(y,z) \in A_{ijk}'$  we have for all $l$ that
$$|(\sum_m \beta_m \partial_{z_m})^q v_{ijk}(y,z)| \leq \delta_{ijkl} |y^{\alpha_{ijk}}| \eqno (3.14)$$
Combining $(3.14)$ with $(3.11)$ gives $(3.9a)-(3.9b)$ and we are done.\hfill $\square$

\noindent {\bf Remark.} By Lemma 3.3 and the fact that the $|z_m|$ are bounded below, in $(3.9a)-(3.9b)$ we can replace
$y^{\alpha_{ijk}}$ or $y^{\alpha_{0jk}}$ by $m \circ \gamma_k(y,z)$ or $m \circ \gamma_k(y)$ respectively, where $m(x) = x^v$
for any vertex $v$ of $N(f)$ on $F_{ij}$. Note the constants $\delta_{ijkl}$ may change.

The final lemma of this section shows how one can form partitions of unity into functions that can be blownup by the
$\gamma_k$'s into smooth functions supported on sets that are approximately subsets of $A_{ijk}'$, and on whose support $(3.9a)$ or $(3.9b)$ holds.

\noindent {\bf Lemma 3.5.} Let $\phi(x) \in C_c(E)$, where $E$ is the cube $\{x \in K_0^n: 0 < |x_l| < C_n^{-1}$ for all $l\}$
which was partitioned into the $A_{ij}$ in section 2. Then except on a set of measure zero one can write $\phi(x) = \sum_{ijkl}\phi_{ijkl}(x)$, where each $\phi_{ijkl}$ is supported in $E$ and each $\rho_{ijkl} = \phi_{ijkl} \circ \gamma_k$ is a smooth function defined on $K_0^n$ which
extends to a smooth function on all of $K^n$ such that the following hold.

\noindent {\bf a)} Suppose $i > 0$. Let $Y = \{y: |y_l| < 2$ for all $l\}$, $Y' = \{y: |y_l| < {1 \over 2}C_n^{-1} $ for all $l\}$, and $Z = \{z: {1 \over 2}C_i^{-1} < |z_m| < 2$ for all $m\}$. For some nonempty open $Z' \subset Z$ we have 
$$ Y' \times Z' \subset supp(\rho_{ijkl}) \subset Y \times Z  \eqno (3.15)$$
Furthermore, for some $0 \leq p \leq o_{ij}$, on $supp(\rho_{ijkl})$ equation $(3.9a)$ holds, possibly with a different constant.

\noindent {\bf b)} Suppose $i = 0$. Let $Y = \{y: |y_l| < 2$ for all $l\}$, $Y' = \{y: |y_l| < {1 \over 2}C_n^{-1} $ for all $l\}$. Then
$$Y' \subset supp(\rho_{0jkl}) \subset Y \eqno (3.16)$$
Furthermore, on $supp(\rho_{0jkl})$ equation $(3.9b)$ holds, possibly with a different constant.

\noindent {\bf Proof.} Let $\psi(x)$ be a smooth nonnegative function on $K$ such that $\psi(x)$ is supported on $|x| > 1 + \epsilon$ such that $\psi(x) = 1$ for $|x| < 1 - \epsilon$, where $\epsilon$ will be dictated by our arguments. In the case where $K$ is a $p$-adic field
we can take $\psi(x)$ to just be the characteristic function of $\{x: |x| \leq 1\}$. We consider the sequence of blowups that 
occurred when the proof of Lemma 3.1 was applied to the vertices of $N(f)$, creating the cubes  $E_k$ that were used in defining the
$A_{ijk}$. Suppose the first blowup was in the $x_p$ and $x_q$ variables. For $\phi(x) \in C_c(E)$, we write $\phi(x) = \phi_1(x) + 
\phi_2(x)$, where $\phi_1(x) = \psi({x_q \over x_p})\phi(x)$ and where $\phi_2(x) =  (1 - \psi({x_q \over x_p}))\phi(x)$. 

Doing the blowup turning $x_q$ into $x_px_q$ while leaving the other variables fixed turns $\phi_1(x)$ into the smooth 
function $\psi(x_q)\phi(x_1,...,x_{q-1}, x_px_q, x_{q+1},...,x_n)$, and similarly doing the blowup converting $x_p$ into $x_px_q$
while leaving the other variables fixed turns $\phi_2(x)$ into the smooth function $(1 - \psi({1 \over x_p})) \phi(x_1,...,x_{p-1}, x_px_q, x_{p+1},...,x_n)$. Note both functions are supported on $B(0,1 + 2\epsilon) = \{x \in K^n: |x_l| < 1 + 2\epsilon$ for all $\l\}$.
We next perform the analogous procedure on the blown-up $\phi_1(x)$ and $\phi_2(x)$ with respect to the second blowup of the resolution process of
Lemma 3.1, obtaining four smooth functions, each supported in the same cube $B(0,1 + 2\epsilon)$. After finitely many steps we are done, and we obtain a finite collection
 $\rho_1(x),...,\rho_q(x)$ of smooth functions such that each $\rho_k(x)$ is supported on a cube slightly larger than $E_k$, and
such that outside a set of measure zero $\phi(x) = \sum_{k=1}^q \rho_k \circ \gamma_k^{-1}(x)$. 

By continuity, $(3.9a)-(3.9b)$ will still hold (possibly with different constants) if $A_{ijk}' =
\gamma_k^{-1}(A_{ij}) \cap E_k $ is enlarged to $A_{ijk}'' = \gamma_k^{-1}(A_{ij}) \cap B(0,1 + 2\epsilon)$ for small enough $\epsilon$, with corresponding slightly larger $B_{ijkl}''$. So we can do a partition of unity and write $\rho_k = \sum_{ijkl}
\rho_{ijkl}'$ such that $(3.9a)$ or $(3.9b)$ holds on the support of $\rho_{ijkl}'$. Furthermore, if $\epsilon$ is small enough,
then $(3.15)$ and $(3.16)$ will hold as a consequence of Lemma 3.2'. (For the left-hand sides, we are using that the $G_k$ contain
$B(0,C_n^{-1})$ since $G_k$ is defined via monomial conditions $|\gamma_{kl}(x)| < C_n^{-1}$).
Thus if we let $\rho_{ijkl}$ be $\rho_{ijkl}'$ times the characteristic function of $K_0^n$, then $\phi(x) = \sum_{ijkl}
\rho_{ijkl} \circ \gamma_k^{-1}(x)$ gives the desired decompositon of $\phi(x)$ and we are done. \hfill $\square$

\noindent {\heading 4. Proof of Theorem 1.1.}

We now commence the proof of Theorem 1.1. The analytic diffeomorphism conditions stipulated by Theorem 1.1 will follow from the Jacobian determinant
condition in part c) because all the maps we will use will be analytic diffeomorphisms on their domains except where the Jacobian determinant 
vanishes. Thus we need not discuss these conditions in the proof.

Let $f(x) = \sum_\alpha f_{\alpha}x^{\alpha}$ be a function represented by a (nontrivial) power series convergent on 
some neighborhood of the origin in $K^n$, with $f(0) = 0$. We will prove Theorem 1.1 by induction on the dimension $n$.
Since the $n=1$ case is trivial, we always assume $n \geq 2$.

Let $\alpha$ be a multiindex with $|\alpha|$ minimal such that $f_{\alpha}$ is nonzero. Let $m = |\alpha|$. As mentioned 
in the proof of Lemma 3.4, the derivatives of order $m$ are spanned by directional derivatives with rational coefficients, 
so we may let $\partial_v = \sum_k \beta_k \partial_{x_k}$, $\beta_k \in \Q$, be a directional derivative such that 
$\partial_v^m f(0) \neq 0$. We may do a linear change of variables with rational coefficients and assume that $v$ is the $x_n$ 
direction. Hence without loss of generality we may assume that
$$\partial_{x_n}^m f(0) \neq 0 \eqno (4.1a)$$
Since $m$ is minimal, for $l < m$ we have
$$\partial_{x_n}^l f(0) = 0 \eqno (4.1b)$$
For our given $n$, we will prove Theorem 1.1  by induction on the $m$ for which $(4.1a)-(4.1b)$ hold. Namely, we assume that either $m = 1$ or that $m > 1$ and 
we know Theorem 1.1 for all $m' < m$. We will actually be inductively proving a slightly stronger statement than Theorem 1.1; in
addition to showing $f \circ \alpha_i(x)$ and the Jacobian of $\alpha_i(x)$ are of the form $c_i(x)m_i(x)$, $c_i(x)$ nonvanishing and $m_i(x)$ a monomial, we will also show that for each coordinate function $g_j(x) = x_j$ for $j < n$ the functions $g_j \circ \alpha_i(x)$ are also of this form.

 By the implicit function theorem applied to $\partial_{x_n}^{m-1} f$ (the proof for the case $K = \R$ works equally well for
arbitrary $K$), there exists some function $g(x_1,...,x_{n-1})$ with a convergent power series on a neighborhood of 
the origin in $K^{n-1}$ such that $g(0) = 0$ and 
$$\partial_{x_n}^{m-1} f(x_1,...,x_{n-1},g(x_1,...,x_{n-1})) = 0 \eqno (4.2)$$
Hence if we let $F(x) = f(x_1,...,x_{n-1},x_n + g(x_1,...,x_{n-1}))$ we have 
$$\partial_{x_n}^{m-1}F(x_1,...,x_{n-1},0) = 0 \eqno (4.3)$$
Although we are not emphasizing it here, in the language of resolution of singularities we are taking a surface resembling a hypersurface of maximal contact and shifting it to the hyperplane $x_n = 0$. Next, we write $F(x)$ as
$$F(x) = h_m(x)x_n^m + \sum_{p < m-1}h_p(x_1,...,x_{n-1})x_n^p \eqno (4.4)$$
There is no $p = m - 1$ term in $(4.4)$ due to the condition $(4.3)$, and $h_m(0) \neq 0$ due to the condition $(4.1a)$. 
We also have $h_p(0) = 0$ for $p < m$ due to $(4.1b)$.
Note that if $m = 1$, the sum in $(4.4)$ is empty, and we have $F(x) = h_1(x)x_n$. Thus we are done with the $m = 1$ case 
to start off the induction.

We now use the induction hypothesis on the dimension $n$ and apply Theorem 1.1 to  $z(x_1,...,x_{n-1}) = 
\prod_{k = 1}^{n-1}x_k  \times \prod_{p < m - 1}h_p(x_1,...,x_{n-1})$. (If an $h_p(x_1,...,x_{n-1})$ is the zero function we exclude it from the product). Then if $\psi(x_1,...,x_{n-1})$ is a cutoff function supported on a sufficiently small neighborhood of the origin of
$K^{n-1}$, $\psi(x_1,...,x_{n-1})$ can be written as a finite sum $\sum_i \psi_i(x_1,...,x_{n-1})$ satisfying the conditions
of Theorem 1.1. So if $\phi(x_1,...,x_n)$ is a cutoff function supported on a sufficiently small neighborhood of the origin of
$K^n$, one can take a $\psi(x_1,...,x_{n-1})$ equal to 1 on a neighborhood of the support of $\phi$, write $\psi = \sum_i \psi_i$
in the above fashion, and then induce a decomposition $\phi = \sum_i \phi_i$ where $\phi_i = \phi \psi_i$. If $\alpha_i$ is one 
of the functions arising from Theorem 1.1 applied to $z(x_1,...,x_{n-1})$ here, then by part b) of Theorem 1.1, $\phi_i \circ \alpha_i(x)= (\phi \circ \alpha_i(x))
( \psi_i \circ \alpha_i(x))$ can be adjusted on the coordinate hyperplanes $\{x_l = 0\}$ for $l < n$ to become a smooth function whose support is 
contained in an open set where  $z \circ \alpha_i (x_1,...,
x_{n-1})$ is of the desired form $c_i(x_1,...,x_{n-1})m_i(x_1,...,x_{n-1})$, where $c_i(x_1,...,x_{n-1})$ is nonvanishing and $m_i(x_1,...,x_{n-1})$ is a monomial. 

It is easy to show that if a product
of functions is of this form, then so is each member of that product. Thus given the definition of $z(x_1,...,x_{n-1})$, the functions
$x_k \circ \alpha_i$ for $k < n$ and each $h_p(x_1,...,x_{n-1})$  is also of the desired form $c_i(x)m_i(x)$ on a neighborhood of the support of $\phi_i(x)$. Correspondingly, switching indices from $i$ to $q$, we write each $h_p \circ \alpha_q (x_1,...,x_{n-1})$ in the form $a_{pq}(x_1,...,x_{n-1})m_{pq}(x_1,...,x_{n-1})$, where $m_{pq}(x_1,...,x_{n-1})$ is either a monomial or the zero function, and each  $a_{pq}(x)$ is nonvanishing on a neighborhood of the support of $\phi_q \circ \alpha_q(x)$. We therefore
can write
$$F \circ \alpha_q(x) =  a_{mq}(x)x_n^m + \sum_{p < m-1}a_{pq}(x_1,...,x_{n-1})m_{pq}(x_1,...,x_{n-1})x_n^p \eqno (4.5)$$
Here $a_{pq}(x)$ also is nonvanishing on a neighborhood of the support of $\phi_q \circ \alpha_q$. In order to prove 
Theorem 1.1, it suffices to show that for each $q$ one can write $\phi_q \circ \alpha_q(x)$ as a sum $\sum_j \eta_j(x)$ of cutoff functions for which
 $F \circ \alpha_q(x)$ and each $x_k \circ \alpha_q(x)$ for $k < n$ satisfies Theorem 1.1 with these $\eta_j(x)$. Since each $x_k \circ \alpha_q(x)$
is of the canonical form $c(x)m(x)$, $m(x)$ a monomial and $c(x)$ nonvanishing, this follows from showing that $F \circ \alpha_q(x)$ and each coordinate function $x_k$ for $k < n$ satisfies Theorem 1.1  with these $\eta_j(x)$, which is what we will prove.

\noindent {\bf Claim.} In order to prove the inductive step of Theorem 1.1 it suffices to prove that for any point $x_0$ in the support of any $\phi_q \circ \alpha_q(x)$ there is an open set $U_{x_0}$ centered at $x_0$ such that Theorem 1.1 holds simultaneously for
$F \circ \alpha_q(x + x_0)$ and each function $x_k$ for $k < n$.

\noindent {\bf Proof.} Suppose we have proven the above. We can use a partition of unity
to write  $\phi_q \circ \alpha_q(x)$ as a sum $\sum_j \eta_j(x)$ of cutoff functions such that each $\eta_j$ is supported in one of
the $U_{x_0}$, so that we may apply Theorem 1.1 to $F \circ \alpha_q(x + x_0)$ and the functions $x_k$ for $k < n$
using the cutoff functions $\eta_i(x + x_0)$. This gives most what we are trying to prove; what remains to do is first proving 
that the composition of the coordinate
functions $x_k$ for $k < n$ with the shift by $x_0$ and then with the coordinate changes from applying the induction hypothesis to 
$F \circ \alpha_q(x + x_0)$ are of the form $c(x)m(x)$, where $m(x)$ is a monomial and $c(x)$ is nonvanishing, and secondly
proving that the determinant of the composition of all the maps in question is of the desired form.

Assuming the $U_{x_0}$ are chosen sufficiently small, if $k$ and $x_0$ are such that 
$(x_0)_k \neq 0$ then $x_k + (x_0)_k$ has no zeroes in $U_{x_0}$, while if $(x_0)_k = 0$ then 
$x_k + (x_0)_k$ is just $x_k$. So if one proves that all the $x_k$ in the shifted coordinates are monomialized by the coordinate
changes on $F \circ \alpha_q(x + x_0)$, the same will be true of the original functions $x_k$ under the composition of the
shift with these coordinate changes. As for the determinant, by induction hypothesis we used in $n-1$ dimensions, the determinant of $\alpha_q$ is of
the form $c_q(x_1,...,x_{n-1})m_q(x_1,...,x_{n-1})$, $m_q$ a monomial and $c_q$ nonvanishing. Similar to above, composing
 this with the shift by $x_0$ results in a function of the same form as long as $U_{x_0}$ is sufficiently small.
 Since are assuming we know Theorem 1.1 for each 
$x_k$ for $k < n$ in the shifted coordinates, and that the determinant of the coordinate changes on $F \circ \alpha_q(x + x_0)$
are of the form $c_q(x)m_q(x)$ for a monomial $m_q$ and a nonvanishing $c_q$, the chain rule for
determinants now implies that the determinant of the composition of all coordinate changes is also of this form $c_q(x)m_q(x)$. Thus the claim is proven. \hfill$\square$

We now proceed to show that for $x_0$ in the support of $\phi_q \circ \alpha_q(x)$,  Theorem 1.1 holds for $F \circ \alpha_q(x + x_0)$ and each coordinate function $x_k$ with $k < n$.
By $(4.5)$, $\partial^{m-1} (F \circ \alpha_q)(x_0) \neq 0$ for $x_0$ in the support of $\phi_q \circ \alpha_q(x)$ except at those $x_0 = ((x_0)_1,...,(x_0)_n)$ for which $(x_0)_n = 0$. So when $(x_0)_n \neq 0$ one can apply the
induction hypothesis and get that Theorem 1.1 holds for $F \circ \alpha_q(x + x_0)$ and each coordinate function $x_k$ with 
$k < n$. Similarly, for $p < m-1$ for which the sum $(4.5)$ has a nonzero term, one has that
$\partial^p (F \circ \alpha_q) (x_0) \neq 0$ unless $a_{pq}((x_0)_1,...,(x_0)_{n-1}) = 0$, and one again can apply 
the induction hypothesis to get that Theorem 1.1 holds for $F \circ \alpha_q(x + x_0)$ and each coordinate function $x_k$ with 
$k < n$.  Thus
we need only consider $x_0$ such that $(x_0)_n = 0$ and each  $a_{pq}((x_0)_1,...,(x_0)_{n-1}) = 0$. In this case, $F \circ \alpha_q(x + x_0)$ is once again of the form $(4.5)$. So without loss of generality, we take $x_0 = 0$ in our subsequent 
arguments and look at $F \circ \alpha_q(x)$.

Let $g(x) = F \circ \alpha_q(x)$. Note that $(4.5)$ implies that $N(g)$ has a special form. Namely, there 
is a vertex at $(0,...,0,m)$, while there is at most one vertex at height $h$ for any $h < m$. Furthermore, if $\sum_{\beta}
g_{\beta}x^{\beta}$ denotes the Taylor expansion of $g(x)$ at the origin, then $g_{\beta} = 0$ for all $\beta = (\beta_1,...,\beta_n)$ with $\beta_n = m-1$. In particular, there is no vertex at height $m-1$. These imply the following important fact:

\noindent {\bf Lemma 4.1.} If $F$ is any compact face of $N(g)$, then any zero of $g_F(x)$ (cf Definition 2.2) in $K_0^n$ has order at most $m-1$.

\noindent {\bf Proof.} Let $F$ be a compact face of $N(g)$. Since there is at most one vertex at any given height, we may let
$v$ be the vertex of $N(g)$ such that the $n$th component $v_n$ is maximal. By the above considerations, either $v_n \leq m - 2$
or $v_n = m$. In the former case $\partial_{x_n}^{v_n}g_F(x)$ is a monomial which therefore doesn't vanish on $K_0^n$.
Hence any zero of $g_F(x)$ in $K_0^n$ has order at most $v_n \leq m - 2$. In the latter case, since $g_{\beta} = 0$ for all $\beta$ with $\beta_n = m-1$, we have that $\partial_{x_n}^{m-1}g_F(x)$ is just $x_n$ and therefore any zero of $g_F(x)$ in
$K_0^n$ has order at most $m - 1$. Hence Lemma 4.1 is proven. \hfill $\square$

We move to the next stage of our argument. Recall we have reduced to showing Theorem 1.1 for the coordinate functions $x_k$ for $k < n$ together with $g(x) = F \circ \alpha_q(x)$ satisfying $(4.5)$.
 We apply Lemma 3.5 to $g(x)$, which tells us that for
$\eta(x)$ supported in a sufficiently small neighborhood of the origin we may write $\eta(x) = \sum_r \eta_r(x)$ where
each $\eta_r(x)$ is one of the functions denoted by $\phi_{ijkl}(x)$ in Lemma 3.5. For each $r$ there is a finite composition
of blowups $\gamma_r(x)$ such that each $\eta_r \circ \gamma_r$ can be adjusted on a set of measure zero to become
 a smooth function. Furthermore, Lemma 3.5 says that on the support of $\eta_r \circ \gamma_r$ one of two possibilities 
occurs. In the first possibility (the $i=0$ case of Lemma 3.5), $g \circ \gamma_r(y)$  can be written in the form $y_1^{\alpha_1}...y_n^{\alpha_n}h(y)$ where by $(3.9b)$ there is a constant $\delta > 0$ such that $|h(y)| >
\delta$.  In the second possibility, corresponding to $i > 0$, we can split into $y$ and $z$ variables such that
 $g \circ \gamma_r(y,z)$  can be written as $y_1^{\alpha_1}...y_p^{\alpha_p}h(y,z)$, where by $(3.15b)$ 
some $z$-directional derivative of $h(y,z)$ of order $\leq m - 1$ is bounded below in magnitude by some $\delta > 0$. Doing a
 linear change of variables in $z$ if necessary, we may always assume it is the $z_n$ derivative, so that $|\partial_{z_n}^a
h(y,z)| > \delta$ on the support of $\eta_r \circ \gamma_r$ for some $0 \leq a \leq m -1$. The condition that $a \leq m -1$ comes from Lemma 4.1 and the fact that in Lemma 3.5 the order of this derivative is at most the maximum order of the zeroes 
of the functions $g_F(x)$ on $K_0^n$.

In the first case, we are already done; $g \circ \gamma_r(y)$ is of the required form, and furthermore since the components
of $\gamma_r(y)$ are all monomials, the composition of each coordinate function $x_k$ with $\gamma_r(y)$ is a monomial, as 
is the Jacobian determinant of $\gamma_r$.
Moving on to the second case, if $a = 0$ then we are done for the same reasons as in the first case. So we assume $a > 0$.
If $j_k(x)$ denotes the coordinate function $x_k$, then each $j_k \circ \gamma_r(y,z)$ is a 
monomial in $y$ and $z$. By $(3.15a)$, the $|z_m|$ are bounded above and below on the support of $\eta_r \circ \gamma_r$,
so on this support we can write $j_k \circ \gamma_r(y,z) = c_{kr}(y,z)m_{kr}(y)$, where $m_{kr}(y)$ is a monomial in the $y$ variables
only and where $c_{kr}(y,z)$ is nonvanishing. 

We now apply the induction hypothesis. Namely, we do a partition of unity and write 
$\eta_r \circ \gamma_r(y,z) = \sum_s \sigma_{rs}(y,z)$, where $\sigma_{rs}(y,z)$ is a cutoff function such that the induction hypothesis applies to $h_1(y,z) = h(y +y_0, z + z_0)$ with the cutoff function $\sigma_{rs}(y + y_0, z + z_0)$. We can do this
since
$h_1(y, z)$ everywhere has a nonvanishing $z_n$ derivative of order at most $m - 1$. Note that if the support of $\sigma_{rs}$ 
is appropriately small, which we may assume, then in the expression $g \circ \gamma_r(y + y_0,
z + z_0)= (y_1 + (y_0)_1)^{\alpha_1}...(y_p + (y_0)_p)^{\alpha_p}h_1(y,z)$ the factors $(y_l + (y_0)_l)^{\alpha_l}$ for $(y_0)_l \neq 0$ are never zero on the support of  $\sigma_{rs}(y + y_0, z + z_0)$, while the for the remaining $l$'s the factor
$(y_l + (y_0)_l)^{\alpha_l}$ is just $y_l^{\alpha_l}$. 

Thus on the support of $\sigma_{rs}(y + y_0, z + z_0)$, the function $g \circ \gamma_r(y + y_0,
z + z_0)$ can be written in the form
$g_1(y,z) = c(y)m(y)h_1(y,z)$ for a monomial $m(y)$ and nonvanishing $c(y)$. We now apply the induction hypothesis to $h_1(y,z)$, simultaneously monomializing $h_1(y,z)$, all $y_l$ variables, and  and all $z_l$ variables other than $l = n$. If $\beta_j(x)$ are the resulting coordinate changes then not only is $h_1 \circ \beta_j(x)$ of the desired form $c(x)m(x)$, but 
 each $y_l \circ \beta_j(x)$ is also of the desired form. Thus $g_1 \circ \beta_j(x)$ is of the form $c(x)m(x)$ as well. Thus $g(x)$
transforms in the desired fashion under the composition of $\gamma_r$, the shift by $(x_0,y_0)$ and $\beta_j$.

Similarly, if $j_k \circ \gamma_r(y,z)$ is as in three paragraphs ago, then assuming
the support of $\sigma_{rs}$ is sufficiently small 
$j_k \circ \gamma_r(y + y_0, z + z_0)$ is of the form $c_{kr}(y,z)m_{kr}(y)$ for a monomial $m_{kr}(y)$ and nonvanishing $c_{kr}(y,z)$ 
because as we saw $j_k \circ \gamma_r(y,z)$ was of this form, although
the new $m_{kr}(y)$ will no longer depend on any $y_l$ variable for which $(y_0)_l$ is nonzero, similar to two paragraphs ago.
Thus since all $y_l$ variables transform into functions of the form $c(x)m(x)$ under the coordinate changes $\beta_j$ by 
induction hypothesis, the 
composition of $j_k \circ \gamma_r(y + y_0, z + z_0)$ with each $\beta_j$ will also transform into functions of this $c(x)m(x)$ form.
Thus we see that the coordinate functions $x_k$ transform into functions of the desired form under the composition of $\gamma_r$,
the shift by $(y_0,z_0)$, and $\beta_j$. 

We have now seen that $g$ and each coordinate function $x_k$ transform in the desired fashion under the composition of $\gamma_r$,
the shift by $(y_0,z_0)$, and $\beta_j$. We still must show that the determinant of the composition of 
the above coordinate changes is of the  desired form $c(x)m(x)$. Note that the determinant of $\gamma_k(y,z)$ is a monomial since its components are monomials. Since the $|z_l|$ are bounded below, this determinant can be written as $c(y,z)m(y)$ for a monomial $m(y)$ and
nonvanishing $c(y,z)$. The composition of $\gamma_k(y,z)$ with the shift by $(y_0, z_0)$ will have the same determinant, which
will also be of the form $c(y,z)m(y)$ in the shifted coordinates if the support of $\sigma_{rs}(y,z)$ is appropriately small, similar
to the last paragraph as well as three paragraphs ago.
By induction hypothesis, each $j_k \circ \beta_j(x)$ is of the desired form $c(x)m(x)$, as is the determinant of
$\beta_j$. Plugging these facts into the chain rule for determinants gives that the composition of all 
coordinate changes also is of this form $c(x)m(x)$.

The proof of Theorem 1.1 is now complete; if $\phi(x)$ is a cutoff function as in the statement of Theorem 1.1, $\tau_i(x)$ for $ i = 1,...,N$ denote all the cutoff functions in the final blown-up coordinates  produced
by the above argument, and $\alpha_i(x)$ denotes the composition of coordinate changes corresponding to $\tau_i(x)$,
then we define $\tau_i'(x)$ to be $\tau_i(x)$ times the characteristic function of $K_0^n$ and 
$\phi_i(x) = \tau_i' \circ \alpha_i^{-1}(x)$.  Then the $\phi_i(x)$ satisfy the conditions of Theorem 1.1 and 
we are done. \hfill $\square$

\noindent {\bf  Example 1.} Suppose after being put in the form $(4.5)$,  $F \circ \alpha_q(x)$ has nondegenerate principal part in the sense of Varchenko [V], meaning that for each face compact face $G$ of $N(F \circ \alpha_q)$ the polynomial $(F \circ \alpha_q)_G(x)$ has zeroes
of order at most 1 on $K_0^n$. Then the rest of the proof of Theorem 1.1 proceeds as follows.

One first does the constructions of section 2, subdividing a neighborhood of the origin into the sets $A_{ij}$ on which the 
terms of the Taylor expansion of $F \circ \alpha_q(x)$ whose exponents lie on the face $F_{ij}$ dominate. Next one 
 does the constructions of section 3 on $F \circ \alpha_q(x)$ and accordingly writes a cutoff function 
$\phi(x)$ in the form $\phi(x) = \sum_{ijkl} \phi_{ijkl}(x)$ as in Lemma 3.5. The support of $\phi_{ijkl} \circ \gamma_k(x)$ will
be contained in a set which is approximately the  blowup of $A_{ij}$ by $\gamma_k$, and by Lemma 3.5, on the support of $\phi_{ijkl} \circ \gamma_k(x)$ the function $F \circ \alpha_q \circ \gamma_k(x)$ can be written as $y_1^{a_1}...y_p^{a_p}F_{ijkl}(y,z)$, where on the support of $\rho_{ijkl}(x)$
either $F_{ijkl}(y,z)$ does not vanish or has nonzero gradient in the $z$ variables wherever $F_{ijkl}(y,z) = 0$.

Now one is
effectively done resolving the singularities of $f(x)$; $F_{ijkl}(y,z)$ is comparable to the monomial $y_1^{a_1}...y_p^{a_p}$
on a neighborhood of any $(y,z)$ where $F_{ijkl}(y,z) \neq 0$, while by the implicit function theorem, on a neighborhood of any zero of $F_{ijkl}(y,z)$ locally one can do a quasitranslation in a $z$ variable to convert $F_{ijkl}(y,z)$ into a function of the form
$y_1^{a_1}...y_p^{a_p}z_i c(y,z)$ where $c(y,z)$ is nonvanishing. Note that the determinant of the composition of the two 
coordinate changes is also comparable to a monomial in this situation; $\gamma_k$ has monomial components and therefore 
its determinant is a monomial, and the
 quasitranslation in $z_n$ will not change this since $|z_n|$ is bounded away from $0$.

Note that if $f(x)$ started out having nondegenerate principal part, then one could jump directly to performing the constructions
of sections 2 and 3 on $f(x)$, and the blown up $f \circ \gamma_k(x)$ will have the above properties. So after doing a single quasitranslation the zero set of $f \circ \gamma_k(x)$ will be resolved.

\noindent {\bf Example 2.} Suppose we are in two dimensions, and $f(x,y)$ satisfies $(4.1a)-(4.1b)$, so that $f(x,y)$ has a zero
of order $m$ at the origin and has nonvanishing $m$th $y$ derivative. After a  quasitranslation in the $y$ variable $f(x,y)$ gets
into the form $(4.4)$. We do not have to do a resolution of singularities in $n-1$ dimensions to get in the form $(4.5)$; in two dimensions any $F(x,y)$ satisfying $(4.4)$ is already in the form $(4.5)$. 

One now does the constructions of Section 3 on $F(x,y)$. Then given any
cutoff function $\phi(x,y)$ supported on a small enough neighborhood of the origin, we can write $\phi(x,y) = \sum_{ijkl} \phi_{ijkl}(x,y)$
 as in Lemma 3.5, and we look at the functions $F \circ \gamma_k(x,y)$ on the support of the $\rho_{ijkl}(x,y) = \phi_{ijkl} \circ 
\gamma_k(x,y)$. If $\rho_{ijkl}(x,y)$ corresponds to an $A_{ijk}'$ with $i = 0$, then by Lemma 3.5 b), on the support of $\rho_{0jkl}(x,y)$ the function $F \circ \gamma_k(x,y)$ has already been monomialized and we are done.  

If $\rho_{ijkl}(x,y)$ corresponds to an $A_{ijk}'$ with $i = 1$, then by Lemma 3.5a), on the support of
$\rho_{1jkl}(x,y)$, the function $F \circ \gamma_k(x,y)$ is of the form $x^p g(x,y)$  and $g(x,y)$ has a nonvanishing $y$ derivative of order at most $m-1$. Furthermore, $\rho_{1jkl}(x,y)$ is supported away from $y = 0$ in this situation. One then proceeds to resolve $g(x,y)$ and $x$ simultaneously, which can be done using the induction hypothesis.

\noindent {\heading 5. References.}

\noindent [AGV] V. Arnold, S Gusein-Zade, A Varchenko, {\it Singularities of differentiable maps
Volume II}, Birkhauser, Basel, 1988.\parskip = 3pt\baselineskip = 3pt

\noindent [BM1] E. Bierstone, P. Milman, {\it Semianalytic and subanalytic sets}, Inst. Hautes
Etudes Sci. Publ. Math. {\bf 67} (1988) 5-42.

\noindent [BM2] E. Bierstone, P. Milman, {\it Canonical desingularization in characteristic zero by blowing up the maximum
 strata of a local invariant.} Invent. Math. {\bf 128} (1997), no. 2, 207-302.

\noindent [CoGrPr] T. Collins, A. Greenleaf, M. Pramanik, {\it A multi-dimensional resolution of singularities with applications to analysis}, to appear, Amer. J. of Math.

\noindent [EH] S. Encinas, H. Hauser, {\it Strong resolution of singularities in characteristic
zero,} Comment. Math. Helv., {\bf 77} (2002), 421-445.

\noindent [EV1] S. Encinas, O. Villamayor, {\it Good points and constructive resolution of 
singularities.} Acta Math. {\bf 181} (1998), no. 1, 109-158. 

\noindent [EV2] S. Encinas, O. Villamayor, {\it A new proof of desingularization over fields of 
characteristic zero.} Proceedings of the International Conference on Algebraic Geometry and 
Singularities (Spanish) (Sevilla, 2001). Rev. Mat. Iberoamericana {\bf 19} (2003), no. 2, 339-353. 

\noindent [G1] M. Greenblatt, {\it A Coordinate-dependent local resolution of singularities and 
applications},  J. Funct. Anal.  {\bf 255}  (2008), no. 8, 1957-1994.

\noindent [G2] M. Greenblatt, {\it Resolution of singularities, asymptotic expansions of oscillatory 
integrals, and related Phenomena}, J. Analyse Math. {\bf 111} (2011) no. 1, 221-245.  

\noindent [G3] M. Greenblatt, {\it Oscillatory integral decay, sublevel set growth, and the Newton
polyhedron}, Math. Annalen {\bf 346} (2010), no. 4, 857-895.

\noindent [G4] M. Greenblatt, {\it Applications of an elementary resolution of singularities algorithm to exponential sums
and congruences modulo $p^n$}, to appear, Israel J. of Math.

\noindent [H1] H. Hironaka, {\it Resolution of singularities of an algebraic variety over a field of characteristic zero I}, 
 Ann. of Math. (2) {\bf 79} (1964), 109-203.

\noindent [H2] H. Hironaka, {\it Resolution of singularities of an algebraic variety over a field of characteristic zero II},  
Ann. of Math. (2) {\bf 79} (1964), 205-326. 

\noindent [J] H. E. W. Jung, {\it Darstellung der Funktionen eines algebraischen K\"orpers zweier unabh\"angiger 
Ver\"anderlichen $x, y$ in der Umgebung einer Stelle $x = a, y = b$}, J. Reine Agnew. Math. {\bf 133} (1908), 289-314.

\noindent [K] J. Kollar,  {\it Lectures on resolution of singularities}, Annals of Mathematics Studies {\bf 166} Princeton University Press, Princeton, NJ, 2007. vi+208 pp. 

\noindent [Pa1] A. Parusinski, {\it Subanalytic functions}, Trans. Amer. Math. Soc. {\bf 344} (1994), 583-595.

\noindent [Pa2] A. Parusinski, {\it On the preparation theorem for subanalytic functions}, New developments in
singularity theory (Cambridge 2000), Nato Sci. Ser. II Math. Phys. Chem., 21, Kluwer Academic
Publisher, Dordrecht (2001), 193-215.

\noindent [PS] D. H. Phong, E. M. Stein, {\it The Newton polyhedron and
oscillatory integral operators}, Acta Mathematica {\bf 179} (1997), 107-152.

\noindent [PSSt] D. H. Phong, E. M. Stein, J. Sturm, {\it Multilinear level
set operators, oscillatory integral operators, and Newton diagrams}, Math.
Annalen, {\bf 319} (2001), 573-596.

\noindent [St] E. Stein, {\it Harmonic analysis; real-variable methods, orthogonality, and oscillatory 
integrals}, Princeton Mathematics Series Vol. 43, Princeton University Press, Princeton, NJ, 1993.

\noindent [Su] H. Sussman, {\it Real analytic desingularization and subanalytic sets:
an elementary approach}, Trans. Amer. Math. Soc. {\bf 317} (1990), no. 2, 417-461.

\noindent [V] A. N. Varchenko, {\it Newton polyhedra and estimates of oscillatory integrals}, Functional 
Anal. Appl. {\bf 18} (1976), no. 3, 175-196.

\noindent [W] J. Wlodarczyk, {\it Simple Hironaka resolution in characteristic zero.} J. Amer. Math. Soc. {\bf 18} (2005), no. 4, 779-822.

\line{}
\line{}

\noindent Department of Mathematics, Statistics, and Computer Science \hfill \break
\noindent University of Illinois at Chicago \hfill \break
\noindent 322 Science and Engineering Offices \hfill \break
\noindent 851 S. Morgan Street \hfill \break
\noindent Chicago, IL 60607-7045 \hfill \break
\noindent greenbla@uic.edu \hfill \break

\end